\documentclass[a4paper]{article}
\usepackage{a4wide}
\usepackage{longtable}
\usepackage{booktabs} 

\usepackage[ruled]{algorithm2e} 

\usepackage{rotate}
\usepackage{amsfonts}
\usepackage{amssymb}
\usepackage[latin1]{inputenc}
\usepackage{epic, eepic}
\usepackage{subfigure}
\usepackage{url}
\usepackage{amsfonts,amsmath}
\usepackage{subfigure}
\usepackage{pstricks,pst-node,psfrag}
\usepackage{graphicx}
\usepackage{hyperref}
\usepackage{url}

\definecolor{green}{rgb}{0.0,0.7,0.0}
\definecolor{lightgreen}{rgb}{0.0,1.0,0.0}
\definecolor{darkgreen}{rgb}{0.0,0.5,0.0}

\definecolor{lightblue}{rgb}{0.4,0.6,1.0}
\definecolor{darkblue}{rgb}{0.0,0.0,0.6}

\newcommand{\nref}[1]{(\ref{#1})}

\newcommand{\R}{{\mathbb{R}}}
\newtheorem{exm}{Example}

\usepackage{ulem}

\newcommand{\radd}[1]{#1}
\newcommand{\rdel}[1]{}
\newcommand{\rsub}[2]{#2}

\begin{document}
\title{High Performance Block Incomplete LU Factorization}


\author{Matthias Bollh\"ofer
\thanks{%
  Institute of Computational Mathematics,  TU Braunschweig, 
  Universit\"atsplatz 2,
    D-38106 Braunschweig,
 Germany,  m.bollhoefer@tu-bs.de}
 \and
Olaf Schenk
\thanks{%
  Institute of Computational Science, Faculty of Informatics, Universit\`a della Svizzera italiana, Lugano, Switzerland, olaf.schenk@usi.ch}
\and Fabio Verbosio
\thanks{%
  Institute of Computational Science, Faculty of Informatics, Universit\`a della Svizzera italiana, Lugano, Switzerland, fabio.verbosio@usi.ch}}

\maketitle

\begin{abstract}
Many application problems that lead to solving linear systems
make use of preconditioned Krylov subspace solvers to compute their solution. 
Among the most popular preconditioning approaches are incomplete
factorization methods either as \rdel{single--level}\radd{single-level} approaches or within
a multilevel framework. 
We will present a block incomplete \rdel{triangular} factorization that is based
on skillfully blocking the system initially and throughout the factorization.
This approach allows for the use of \rdel{cache--optimized}\radd{cache-optimized} dense matrix
kernels such \rdel{level--3}\radd{as level-3} BLAS or LAPACK. We will demonstrate how this
block approach outperforms the scalar method often by orders of magnitude
on modern architectures, paving the way for its prospective use inside various
multilevel incomplete factorization approaches or other applications where the
core part relies on an incomplete factorization.
\end{abstract}

%
%

%
%


\section{Introduction}\label{sect:intro}
Many application problems lead to solving linear systems of type
\[
Ax=b,
\]
where $A$ is an $n\times n$ nonsingular real or complex system \radd{matrix} and $b$ 
is the associated right-hand side. In particular we are interested in the
case where $A$ is \rsub{large--scale}{large-scale} and sparse. The generic way of solving 
these systems nowadays consists of using \rdel{state--of--the--art}\radd{state-of-the-art} sparse direct solvers 
(cf.\rdel{ e.g.}\radd{, e.g.,}~\cite{AmeDKE01,sg:04-fgcs,Dav04,li05}). 
Although high performance sparse direct solvers
a\radd{re} very efficient in many cases, several structured problems\radd{, i.e., problems presenting specific, noticeable sparsity structures,} cause the
direct solver to produce a significant amount of fill-in during the
factorization\radd{,} leading to high memory requirements which \rsub{are often not within 
the scope of the underlying machine}{can exceed the hardware capability}. If these kind of problems can be solved
efficiently, then one has to rely on \rsub{out--of--core}{out-of-core} techniques.
\radd{These techniques rely on memory locations external to the computer's working memory, i.e., disks, in order to overcome hardware limitations; see, e.g.,~\cite{AmeDLR12}.}
\rsub{Besides, the computation time usually increases drastically. In such 
situations
there is a realistic hope that approximate factorization methods in 
combination with Krylov subspace methods offer a more efficient and 
attractive alternative.}{The presence of high fill-in might lead to prohibitive execution time, suggesting the use of approximate factorization strategies in combination with Krylov subspace methods as a valid alternative approach.}
\rsub{Among the most popular 
approximate factorization methods
are in particular those based on incomplete $LU$ factorization with
several variations~[Saad 2003b] and more recent, nowadays 
$ILU$ approaches inside a multilevel framework (see e.g.~[Bollhöfer and Saad 2006; Saad and Suchomel 2002; Xi et al. 2016]
just to mention some of them).}{Among the most popular approximate factorization methods, we mention those based on the incomplete $LU$ factorization~\cite{Saa03} and the more recently developed $ILU$ approaches in multilevel frameworks, such as~\cite{SaaS02,BolS06,XiLS16}.}
For direct $LU$ factorization methods\radd{,} block structured algorithms\radd{,}
such as multifrontal methods or those based on supernodes\radd{,} have demonstrated
their superiority on modern hardware architectures mainly due to the
usage of dense linear algebra kernels such as \rsub{level--3 BLAS for matrix--matrix}{level-3 BLAS for matrix-matrix}
operations or LAPACK for certain factorization templates. 
Part of the success of direct solvers is obtained from the symbolic
analysis using the (column) elimination tree which is able to predict
dense blocks in advance and to set up the data structures appropriately.
For incomplete $LU$ factorization methods this is usually not possible
except for very few approaches such as \rsub{level--of--fill}{level-of-fill} approaches~\cite{HenRR08}. 
For the symmetric positive definite case, in~\cite{NgPR99} a block incomplete
Cholesky decomposition is computed which uses
a supernodal structure breaking up the 
supernodes into smaller blocks in order to allow refined dropping.
A generic approach to block preconditioning method is introduced
in~\cite{ChoH98}, where a C++ framework is provided
offering \rsub{block--oriented}{block-oriented} preconditioning methods for block structures
defined by the user\rsub{,}{;} one of these is a block\rsub{--}{ }tridiagonal ILU.
%
A very efficient and successful incomplete
Cholesky factorization method was presented in~\cite{GupG10}, where several
aspects\radd{,} such as blocking using the elimination tree or efficient implementation
using dense matrix kernels\radd{,} were put together to eventually end up
in a very robust sparse block incomplete Cholesky factorization method.
Furthermore, a supernodal block incomplete factorization 
approach was presented in~\cite{LiS11}.

In the present paper, our block ILU approach uses several
known components, combines them but also introduces further strategies to 
construct efficient block structures \radd{with blocks of variable size}
for a block ILU factorization method.
Furthermore we improve the block partitioning during the factorization. 
It is the combination of several ingredients
that eventually improves the block ILU method significantly 
over its scalar \rsub{counter part}{counterpart} in many practical applications on modern
computer architectures. Our approach thus generalizes the scalar ILU approach
to a block approach, yet further prospective applications of this approach 
are subject to future research such as using block ILUs within a multilevel
framework.

The paper is organized as follows.
We will briefly review established incomplete factorization methods
(\rsub{S}{s}ection~\ref{sect:problem}) with special focus \rsub{of}{on} the \rsub{so--called Crout--type}{so-called Crout-type}
ILU which is sometimes also referred to \rsub{left--looking}{as left-looking} 
ILU (at least with respect to $L$).
We will demonstrate that this approach
\rsub{is easily upgraded}{can easily be extended} to a block ILU and focus on the major challenges when
switching to a block method.
Section~\ref{sect:preprocessing} is devoted
to provid\rsub{e}{ing} the block structures required to make the block ILU approach
efficient. It comprises techniques to improve diagonal dominance, reduction
of \rsub{fill--in}{fill-in} as well as a priori \radd{variable} block partitioning and aggregating blocks
during the factorization.
Finally we will demonstrate in \rsub{S}{s}ection~\ref{sect:exp} that the combination
of these technologies ends up in a very efficient high performance
incomplete factorization approach which can easily outperform the
traditional ILU by orders of magnitude on modern computer\radd{s} using dense
matrix kernels.

\section{Incomplete Factorization Methods}\label{sect:problem}
The design of preconditioning methods based on incomplete $LU$ factorization
typically relies on efficiently computing approximate triangular factors
without having too much symbolic information \rsub{at}{on} hand. \rsub{Certainly f}{F}or
level-of-fill ILUs one can certainly use information from the elimination
tree~\rdel{(}\cite{HenRR08}\rdel{)}. In contrast to that, \rsub{threshold--based}{threshold-based} ILUs
are hardly able to use this kind of information. Instead, efficiency requires \radd{us}
to either compute \rsub{extremely sparse and robust factors}{significantly sparser factors which remain robust in spite of dropping} or to 
heuristically introduce block structures to increase performance~\cite{CarLS14,GupG10}.
The \rsub{classical}{general} incomplete $LU$ factorization approaches distinguish
how the portions of $L$ and $U$ are to be computed,
e.g.\radd{,} row\rdel{--}wise (also referred to as IKJ variant or known as \radd{the} ILUT~\cite{Saa03})
as one example.

\subsection{The Crout ILU}\label{sect:iluc}
A particularly attractive incomplete factorization approach is the so--called Crout ILU~\cite{EisGSS82,JonP95,LinM99,Saa03,LiSC03} 
 since
it computes the columns of $L$ and rows of $U$ simultaneously only
using the already computed parts of $L$ and $U$. We highlight this
version since we are going to establish our block ILU based on this variant.
Algorithm~\ref{iluc} gives a rough sketch of this variant omitting several
technical details.

\begin{algorithm}
\caption{Crout ILU\radd{.}}\label{iluc}
  \KwIn{ $A\in\R^{n,n}$, drop tolerance $1>\tau>0$. }
  \KwOut{ approximate factors $L,U$. }
    \For{$k=1,2,\dots,n$}{
        $l_{ik}\leftarrow a_{ik}$, for all $i\geqslant k$ and $a_{ik}\not=0$ \;
        \For{$j=1,2,\dots,k-1$ such that $u_{jk}\not=0$}{
            $l_{ik}\leftarrow l_{ik}-l_{ij}u_{jk}$ for all $i\geqslant k$ and $l_{ij}\not=0$ \;
        }
        drop $l_{ik}$ whenever $|l_{ik}|\leqslant \tau |l_{kk}|$, for all $i\geqslant k$ and $l_{ik}\not=0$ \;
        $l_{ik}\leftarrow l_{ik}/l_{kk}$ for all $i\geqslant k$ and $l_{ik}\not=0$ \;
        $u_{ki}\leftarrow a_{ki}$, for all $i\geqslant k$ and $a_{ki}\not=0$ \;
        \For{$j=1,2,\dots,k-1$ such that $l_{kj}\not=0$}{
            $u_{ki}\leftarrow u_{ki}-l_{kj}u_{ji}$ for all $i\geqslant k$ and $u_{ji}\not=0$ \;
        }
        drop $u_{ki}$ whenever $|u_{ki}|\leqslant \tau |u_{kk}|$, for all $i\geqslant k$ and $u_{ki}\not=0$ \;
    } 
\end{algorithm}

An efficient realization of Algorithm~\ref{iluc} \rsub{does certainly require to efficiently}{certainly does require us to} deal with the updates of column $k$ of $L$ (resp\rsub{ectively}{.,} row $k$ of $U$). 
This is usually realized using a \rsub{dense buffer}{auxiliary vector} and two associated index
arrays, whereas the final sparsified row/column $k$ is stored only
in compressed format~\cite{Saa03}. While this is more or less standard,
the more difficult aspect in Algorithm~\ref{iluc} is to access $L$ row\rdel{-}wise
although it \radd{is} stored in compressed sparse column format (similar arguments apply
to $U$). A very elegant way to achieve this is to use additional auxiliary
\rsub{$n$--dimensional}{$n$-dimensional} vectors 
\texttt{L\_head},\texttt{L\_list},\texttt{L\_first} for $L$
and \texttt{U\_head},\texttt{U\_list},\texttt{U\_first} for $U$ \radd{which goes back to \cite{EisGSS82}
and can also be found in \cite{JonP95}}.
A detailed description \radd{of} how these additional vectors have to be used
can be found in~\cite{JonP95,LiSC03}. 
The same kind of data structures can furthermore be used to access the initial
matrix $A$ by columns and by rows simultaneously \radd{which is often used in sparse matrix-matrix multiplication codes}. The only constraint to
make this approach work is to ensure that the nonzero entries in each 
column of $L$ are stored \rsub{in increasing order}{keeping increasing row indices} (similar requirements are
necessary for $U$ and $A$).
In total, if implemented efficiently, the Crout\rdel{--}\radd{ }ILU is an extremely
\rsub{efficient}{effective} incomplete factorization approach, since on one hand
it computes the columns of $L$ and rows of $U$ simultaneously and on the other
hand it is extremely memory\rdel{--}\radd{ }efficient as there is only \radd{a} constant number
of additional auxiliary arrays of length $n$ required in addition to
the factors $L$, $U$ to be computed. In the next section we will describe
how this approach can be easily turned into a block ILU.


\subsection{Crout-type Block ILU}\label{sect:biluc}
\begin{figure}
\begin{center}
\unitlength 1cm
{
\begin{picture}(11.5,8.5)
%
\put(1.02,2.52){\textcolor{green}{\line(0,1){3.46}}}
\put(1.25,2.52){\textcolor{green}{\line(0,1){3.46}}}
\put(1.02,2.52){\textcolor{green}{\line(1,0){0.23}}}
\put(1.02,5.98){\textcolor{green}{\line(1,0){0.23}}}
\put(1.29,2.52){\textcolor{green}{\line(0,1){2.96}}}
\put(1.98,2.52){\textcolor{green}{\line(0,1){2.96}}}
\put(1.29,2.52){\textcolor{green}{\line(1,0){0.69}}}
\put(1.29,5.48){\textcolor{green}{\line(1,0){0.69}}}
\put(2.02,2.52){\textcolor{green}{\line(0,1){2.46}}}
\put(2.48,2.52){\textcolor{green}{\line(0,1){2.46}}}
\put(2.02,2.52){\textcolor{green}{\line(1,0){0.46}}}
\put(2.02,4.98){\textcolor{green}{\line(1,0){0.46}}}
\put(2.52,2.52){\textcolor{green}{\line(0,1){1.96}}}
\put(2.98,2.52){\textcolor{green}{\line(0,1){1.96}}}
\put(2.52,2.52){\textcolor{green}{\line(1,0){0.46}}}
\put(2.52,4.48){\textcolor{green}{\line(1,0){0.46}}}
%
\put(1.02,7.27){\textcolor{black}{\rule{0.25cm}{0.25cm}}}
\put(1.27,6.52){\textcolor{black}{\rule{0.75cm}{0.75cm}}}
\put(2.02,6.02){\textcolor{black}{\rule{0.5cm}{0.5cm}}}
\put(2.52,5.52){\textcolor{black}{\rule{0.5cm}{0.5cm}}}
%
\put(2.52,7.48){\textcolor{blue}{\line(1,0){3.46}}}
\put(2.52,7.27){\textcolor{blue}{\line(1,0){3.46}}}
\put(2.52,7.27){\textcolor{blue}{\line(0,1){0.22}}}
\put(5.98,7.27){\textcolor{blue}{\line(0,1){0.22}}}
\put(3.02,7.23){\textcolor{blue}{\line(1,0){2.96}}}
\put(3.02,6.52){\textcolor{blue}{\line(1,0){2.96}}}
\put(3.02,6.52){\textcolor{blue}{\line(0,1){0.70}}}
\put(5.98,6.52){\textcolor{blue}{\line(0,1){0.70}}}
\put(3.52,6.48){\textcolor{blue}{\line(1,0){2.46}}}
\put(3.52,6.02){\textcolor{blue}{\line(1,0){2.46}}}
\put(3.52,6.02){\textcolor{blue}{\line(0,1){0.46}}}
\put(5.98,6.02){\textcolor{blue}{\line(0,1){0.46}}}
\put(4.02,5.98){\textcolor{blue}{\line(1,0){1.96}}}
\put(4.02,5.52){\textcolor{blue}{\line(1,0){1.96}}}
\put(4.02,5.52){\textcolor{blue}{\line(0,1){0.46}}}
\put(5.98,5.52){\textcolor{blue}{\line(0,1){0.46}}}
%
\put(3.1,3.5){$L$}
\put(5.0,5.1){$U$}
\put(3.1,5.1){$D$}
%
\put(1.02,5.82){\textcolor{green}{\rule{0.23cm}{0.16cm}}}
\put(1.02,5.32){\textcolor{green}{\rule{0.23cm}{0.24cm}}}
\put(1.02,4.42){\textcolor{green}{\rule{0.23cm}{0.08cm}}}
\put(1.02,4.22){\textcolor{green}{\rule{0.23cm}{0.08cm}}}
\put(1.02,3.52){\textcolor{green}{\rule{0.23cm}{0.16cm}}}
\put(1.02,3.32){\textcolor{green}{\rule{0.23cm}{0.08cm}}}
\put(1.29,5.12){\textcolor{green}{\rule{0.69cm}{0.36cm}}}
\put(1.29,4.32){\textcolor{green}{\rule{0.69cm}{0.08cm}}}
\put(1.29,4.22){\textcolor{green}{\rule{0.69cm}{0.16cm}}}
\put(1.29,4.02){\textcolor{green}{\rule{0.69cm}{0.08cm}}}
\put(1.29,3.52){\textcolor{green}{\rule{0.69cm}{0.16cm}}}
\put(1.29,2.82){\textcolor{green}{\rule{0.69cm}{0.08cm}}}
\put(2.02,4.72){\textcolor{green}{\rule{0.46cm}{0.16cm}}}
\put(2.02,3.92){\textcolor{green}{\rule{0.46cm}{0.08cm}}}
\put(2.02,2.62){\textcolor{green}{\rule{0.46cm}{0.08cm}}}
\put(2.02,2.52){\textcolor{green}{\rule{0.46cm}{0.24cm}}}
\put(2.52,4.32){\textcolor{green}{\rule{0.46cm}{0.16cm}}}
\put(2.52,4.02){\textcolor{green}{\rule{0.46cm}{0.08cm}}}
\put(2.52,3.52){\textcolor{green}{\rule{0.46cm}{0.08cm}}}
\put(2.52,2.62){\textcolor{green}{\rule{0.46cm}{0.08cm}}}
\put(2.52,2.52){\textcolor{green}{\rule{0.46cm}{0.16cm}}}
%
\put(2.52,7.27){\textcolor{blue}{\rule{0.16cm}{0.22cm}}}
\put(2.94,7.27){\textcolor{blue}{\rule{0.24cm}{0.22cm}}}
\put(3.80,7.27){\textcolor{blue}{\rule{0.08cm}{0.22cm}}}
\put(4.20,7.27){\textcolor{blue}{\rule{0.08cm}{0.22cm}}}
\put(4.80,7.27){\textcolor{blue}{\rule{0.16cm}{0.22cm}}}
\put(5.18,7.27){\textcolor{blue}{\rule{0.08cm}{0.22cm}}}
%
\put(3.26,6.52){\textcolor{blue}{\rule{0.24cm}{0.70cm}}}
\put(3.92,6.52){\textcolor{blue}{\rule{0.08cm}{0.70cm}}}
\put(4.14,6.52){\textcolor{blue}{\rule{0.16cm}{0.70cm}}}
\put(4.42,6.52){\textcolor{blue}{\rule{0.08cm}{0.70cm}}}
\put(4.92,6.52){\textcolor{blue}{\rule{0.16cm}{0.70cm}}}
\put(5.62,6.52){\textcolor{blue}{\rule{0.08cm}{0.70cm}}}
%
\put(3.88,6.02){\textcolor{blue}{\rule{0.16cm}{0.46cm}}}
\put(4.48,6.02){\textcolor{blue}{\rule{0.08cm}{0.46cm}}}
\put(5.48,6.02){\textcolor{blue}{\rule{0.08cm}{0.46cm}}}
\put(5.72,6.02){\textcolor{blue}{\rule{0.24cm}{0.46cm}}}
%
\put(4.02,5.52){\textcolor{blue}{\rule{0.16cm}{0.46cm}}}
\put(4.40,5.52){\textcolor{blue}{\rule{0.08cm}{0.46cm}}}
\put(4.90,5.52){\textcolor{blue}{\rule{0.08cm}{0.46cm}}}
\put(5.40,5.52){\textcolor{blue}{\rule{0.08cm}{0.46cm}}}
\put(5.82,5.52){\textcolor{blue}{\rule{0.16cm}{0.46cm}}}
\put(7.0,6.9){$\downarrow$}
\put(7.0,6.4){$\downarrow$}
\put(7.0,5.9){$\downarrow$}
\put(6.2,6.4){~~block~\quad~structure}
\put(4.0,7.9){$\rightarrow$}
\put(4.5,7.9){$\rightarrow$}
\put(5.0,7.9){$\rightarrow$}
\put(4.2,8.1){scalar}
\put(4.0,7.6){structure}
\put( 0.3,5.0){$\downarrow$}
\put( 0.0,4.5){scalar}
\put( 0.3,4.0){$\downarrow$}
\put(-0.2,3.5){structure}
\put( 0.3,3.0){$\downarrow$}
\put(1.3,1.9){$\rightarrow$}
\put(1.8,1.9){$\rightarrow$}
\put(2.3,1.9){$\rightarrow$}
\put(1.6,2.1){block}
\put(1.3,1.6){structure}
\put( 9.0,8.1){block rows stored}
\put( 9.0,7.7){as dense matrices}
\put( 9.0,7.37){\vector(1,0){2}}
\put( 9.0,6.86){\vector(1,0){2}}
\put( 9.0,6.22){\vector(1,0){2}}
\put( 9.0,5.72){\vector(1,0){2}}
\put(11.0,7.27){\textcolor{blue}{\rule{0.8cm}{0.22cm}}}
\put(11.0,6.52){\textcolor{blue}{\rule{0.8cm}{0.70cm}}}
\put(11.0,6.02){\textcolor{blue}{\rule{0.56cm}{0.46cm}}}
\put(11.0,5.52){\textcolor{blue}{\rule{0.56cm}{0.46cm}}}
\put(1.05,1.0){$\downarrow$}
\put(1.55,1.0){$\downarrow$}
\put(2.17,1.0){$\downarrow$}
\put(2.67,1.0){$\downarrow$}
\put(1.02,0.0){\textcolor{green}{\rule{0.23cm}{0.8cm} }}
\put(1.29,0.0){\textcolor{green}{\rule{0.69cm}{0.92cm} }}
\put(2.02,0.0){\textcolor{green}{\rule{0.46cm}{0.44cm}}}
\put(2.52,0.0){\textcolor{green}{\rule{0.46cm}{0.56cm}}}
\put( 3.3,0.6){block columns stored}
\put( 3.5,0.2){as dense matrices}
\end{picture}
}
\end{center}
\caption{Sketch of the block structures of $L,D,U$ inside the block ILU.}
\label{ldu-mem}
\end{figure}
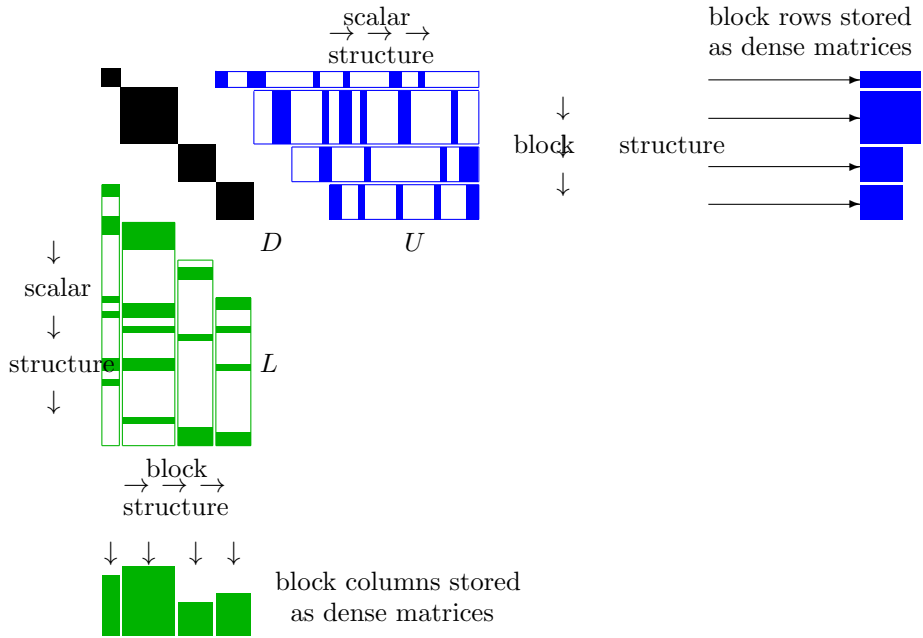

As a first step towards a \rsub{block--structured}{block-structured} ILU we like to point out that
Algorithm~\ref{iluc} can \rdel{be} almost \radd{be} implemented straightforward\radd{ly} in the
same way if the scalar entries are replaced by blocks. 
Formally
only minor changes such as $\|l_{ik}l_{kk}^{-1}\|\leqslant \tau$,
$\|u_{kk}^{-1}u_{ki}\|\leqslant \tau$\radd{,} and $l_{ik}\leftarrow l_{ik}l_{kk}^{-1}$
are necessary for the block version.
In what follows we describe in more detail how our \rdel{block--structured}\radd{block-structured}
version of the Crout ILU is going to be realized.
We assume that the initial matrix is just as regular \radd{a} sparse matrix,
usually without any specific block structure. We may assume that using
some permutation strategy we end up with \radd{a} sparse matrix where at least
a block partitioning for the diagonal blocks is obtained. We will
later comment in more detail about initial preprocessing steps of this
kind. If the \rsub{size of the diagonal blocks}{variable size of each diagonal block}  is given in advance then this
easily imposes a block structure for the block columns of $L$ as well as
the block rows of $U$. However, since we are going to drop entries of
small size, we will maintain a scalar structure for the rows of $L$ and
the columns of $U$. In addition, we will store the dense diagonal blocks
separately in a block diagonal matrix $D$.
This gives a hybrid structure with blocks in one
direction and a scalar representation in the other direction.
The structure is illustrated in Figure~\ref{ldu-mem}.

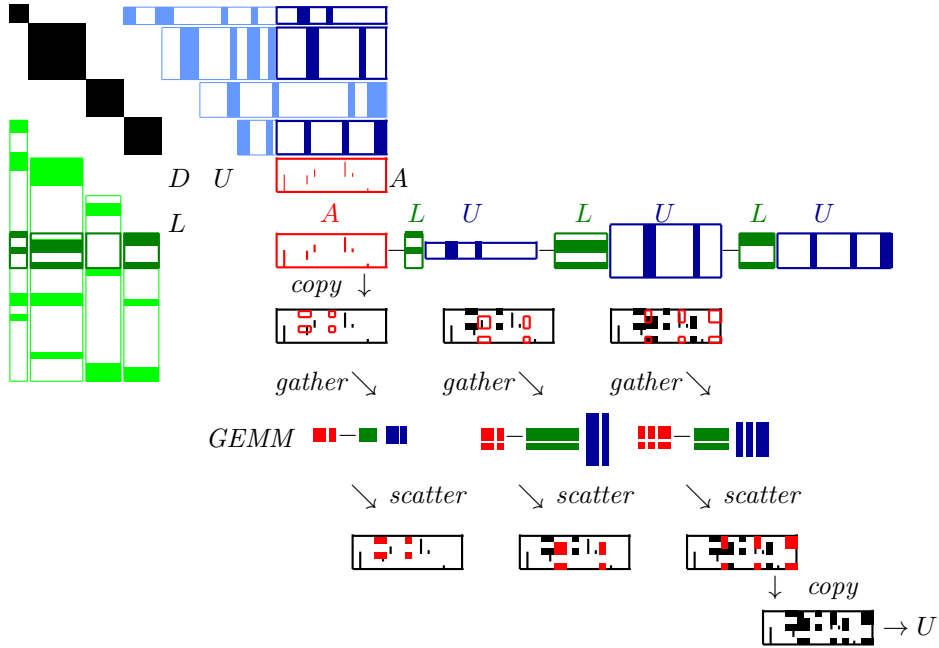
\begin{figure}
\begin{center}
\unitlength 1cm
{
\begin{picture}(11.5,8.5)
%
\put(0.02,3.52){\textcolor{lightgreen}{\line(0,1){3.46}}}
\put(0.25,3.52){\textcolor{lightgreen}{\line(0,1){3.46}}}
\put(0.02,3.52){\textcolor{lightgreen}{\line(1,0){0.23}}}
\put(0.02,6.98){\textcolor{lightgreen}{\line(1,0){0.23}}}
\put(0.29,3.52){\textcolor{lightgreen}{\line(0,1){2.96}}}
\put(0.98,3.52){\textcolor{lightgreen}{\line(0,1){2.96}}}
\put(0.29,3.52){\textcolor{lightgreen}{\line(1,0){0.69}}}
\put(0.29,6.48){\textcolor{lightgreen}{\line(1,0){0.69}}}
\put(1.02,3.52){\textcolor{lightgreen}{\line(0,1){2.46}}}
\put(1.48,3.52){\textcolor{lightgreen}{\line(0,1){2.46}}}
\put(1.02,3.52){\textcolor{lightgreen}{\line(1,0){0.46}}}
\put(1.02,5.98){\textcolor{lightgreen}{\line(1,0){0.46}}}
\put(1.52,3.52){\textcolor{lightgreen}{\line(0,1){1.96}}}
\put(1.98,3.52){\textcolor{lightgreen}{\line(0,1){1.96}}}
\put(1.52,3.52){\textcolor{lightgreen}{\line(1,0){0.46}}}
\put(1.52,5.48){\textcolor{lightgreen}{\line(1,0){0.46}}}
%
\put(0.02,8.27){\textcolor{black}{\rule{0.25cm}{0.25cm}}}
\put(0.27,7.52){\textcolor{black}{\rule{0.75cm}{0.75cm}}}
\put(1.02,7.02){\textcolor{black}{\rule{0.5cm}{0.5cm}}}
\put(1.52,6.52){\textcolor{black}{\rule{0.5cm}{0.5cm}}}
%
\put(1.52,8.48){\textcolor{lightblue}{\line(1,0){3.46}}}
\put(1.52,8.25){\textcolor{lightblue}{\line(1,0){3.46}}}
\put(1.52,8.25){\textcolor{lightblue}{\line(0,1){0.23}}}
\put(4.98,8.25){\textcolor{lightblue}{\line(0,1){0.23}}}
\put(2.02,8.21){\textcolor{lightblue}{\line(1,0){2.96}}}
\put(2.02,7.52){\textcolor{lightblue}{\line(1,0){2.96}}}
\put(2.02,7.52){\textcolor{lightblue}{\line(0,1){0.69}}}
\put(2.98,7.52){\textcolor{lightblue}{\line(0,1){0.69}}}
\put(2.52,7.48){\textcolor{lightblue}{\line(1,0){2.46}}}
\put(2.52,7.02){\textcolor{lightblue}{\line(1,0){2.46}}}
\put(2.52,7.02){\textcolor{lightblue}{\line(0,1){0.46}}}
\put(4.98,7.02){\textcolor{lightblue}{\line(0,1){0.46}}}
\put(3.02,6.98){\textcolor{lightblue}{\line(1,0){1.96}}}
\put(3.02,6.52){\textcolor{lightblue}{\line(1,0){1.96}}}
\put(3.02,6.52){\textcolor{lightblue}{\line(0,1){0.46}}}
\put(4.98,6.52){\textcolor{lightblue}{\line(0,1){0.46}}}
%
\put(2.1,5.5){$L$}
\put(2.7,6.1){$U$}
\put(2.1,6.1){$D$}
%
\put(0.02,6.82){\textcolor{lightgreen}{\rule{0.23cm}{0.16cm}}}
\put(0.02,6.32){\textcolor{lightgreen}{\rule{0.23cm}{0.24cm}}}
\put(0.02,5.42){\textcolor{lightgreen}{\rule{0.23cm}{0.08cm}}}
\put(0.02,5.22){\textcolor{lightgreen}{\rule{0.23cm}{0.08cm}}}
\put(0.02,4.52){\textcolor{lightgreen}{\rule{0.23cm}{0.16cm}}}
\put(0.02,4.32){\textcolor{lightgreen}{\rule{0.23cm}{0.08cm}}}
\put(0.29,6.12){\textcolor{lightgreen}{\rule{0.69cm}{0.36cm}}}
\put(0.29,5.32){\textcolor{lightgreen}{\rule{0.69cm}{0.08cm}}}
\put(0.29,5.22){\textcolor{lightgreen}{\rule{0.69cm}{0.16cm}}}
\put(0.29,5.02){\textcolor{lightgreen}{\rule{0.69cm}{0.08cm}}}
\put(0.29,4.52){\textcolor{lightgreen}{\rule{0.69cm}{0.16cm}}}
\put(0.29,3.82){\textcolor{lightgreen}{\rule{0.69cm}{0.08cm}}}
\put(1.02,5.72){\textcolor{lightgreen}{\rule{0.46cm}{0.16cm}}}
\put(1.02,4.92){\textcolor{lightgreen}{\rule{0.46cm}{0.08cm}}}
\put(1.02,3.62){\textcolor{lightgreen}{\rule{0.46cm}{0.08cm}}}
\put(1.02,3.52){\textcolor{lightgreen}{\rule{0.46cm}{0.24cm}}}
\put(1.52,5.32){\textcolor{lightgreen}{\rule{0.46cm}{0.16cm}}}
\put(1.52,5.02){\textcolor{lightgreen}{\rule{0.46cm}{0.08cm}}}
\put(1.52,4.52){\textcolor{lightgreen}{\rule{0.46cm}{0.08cm}}}
\put(1.52,3.62){\textcolor{lightgreen}{\rule{0.46cm}{0.08cm}}}
\put(1.52,3.52){\textcolor{lightgreen}{\rule{0.46cm}{0.16cm}}}
%
\put(1.52,8.25){\textcolor{lightblue}{\rule{0.16cm}{0.23cm}}}
\put(1.94,8.25){\textcolor{lightblue}{\rule{0.24cm}{0.23cm}}}
\put(2.80,8.25){\textcolor{lightblue}{\rule{0.08cm}{0.23cm}}}
\put(3.20,8.25){\textcolor{lightblue}{\rule{0.08cm}{0.23cm}}}
\put(3.80,8.25){\textcolor{lightblue}{\rule{0.16cm}{0.23cm}}}
\put(4.18,8.25){\textcolor{lightblue}{\rule{0.08cm}{0.23cm}}}
%
\put(2.26,7.52){\textcolor{lightblue}{\rule{0.24cm}{0.69cm}}}
\put(2.92,7.52){\textcolor{lightblue}{\rule{0.08cm}{0.69cm}}}
\put(3.14,7.52){\textcolor{lightblue}{\rule{0.16cm}{0.69cm}}}
\put(3.42,7.52){\textcolor{lightblue}{\rule{0.08cm}{0.69cm}}}
\put(3.92,7.52){\textcolor{lightblue}{\rule{0.16cm}{0.69cm}}}
\put(4.62,7.52){\textcolor{lightblue}{\rule{0.08cm}{0.69cm}}}
%
\put(2.88,7.02){\textcolor{lightblue}{\rule{0.16cm}{0.46cm}}}
\put(3.48,7.02){\textcolor{lightblue}{\rule{0.08cm}{0.46cm}}}
\put(4.48,7.02){\textcolor{lightblue}{\rule{0.08cm}{0.46cm}}}
\put(4.72,7.02){\textcolor{lightblue}{\rule{0.24cm}{0.46cm}}}
%
\put(3.02,6.52){\textcolor{lightblue}{\rule{0.16cm}{0.46cm}}}
\put(3.40,6.52){\textcolor{lightblue}{\rule{0.08cm}{0.46cm}}}
\put(3.90,6.52){\textcolor{lightblue}{\rule{0.08cm}{0.46cm}}}
\put(4.40,6.52){\textcolor{lightblue}{\rule{0.08cm}{0.46cm}}}
\put(4.82,6.52){\textcolor{lightblue}{\rule{0.16cm}{0.46cm}}}
\thicklines
\put(3.53,6.47){\textcolor{red}{\line(1,0){1.46}}}
\put(3.53,6.03){\textcolor{red}{\line(1,0){1.46}}}
\put(3.53,6.03){\textcolor{red}{\line(0,1){0.46}}}
\put(4.97,6.03){\textcolor{red}{\line(0,1){0.46}}}
\put(5.0, 6.10){{$A$}}
\put(0.01,5.48){\textcolor{darkgreen}{\line(1,0){.24}}}
\put(0.01,5.02){\textcolor{darkgreen}{\line(1,0){.24}}}
\put(0.285,5.48){\textcolor{darkgreen}{\line(1,0){0.70}}}
\put(0.285,5.02){\textcolor{darkgreen}{\line(1,0){0.70}}}
\put(1.015,5.48){\textcolor{darkgreen}{\line(1,0){0.47}}}
\put(1.015,5.02){\textcolor{darkgreen}{\line(1,0){0.47}}}
\put(1.515,5.48){\textcolor{darkgreen}{\line(1,0){0.47}}}
\put(1.515,5.02){\textcolor{darkgreen}{\line(1,0){0.47}}}
\put(0.02,5.02){\textcolor{darkgreen}{\line(0,1){0.46}}}
\put(1.98,5.02){\textcolor{darkgreen}{\line(0,1){0.46}}}
\thinlines
\put(3.63,6.03){\textcolor{red}{\line(0,1){0.23}}}
\put(3.93,6.13){\textcolor{red}{\line(0,1){0.12}}}
\put(4.03,6.23){\textcolor{red}{\line(0,1){0.10}}}
\put(4.43,6.23){\textcolor{red}{\line(0,1){0.20}}}
\put(4.53,6.23){\textcolor{red}{\line(0,1){0.05}}}
\put(4.73,6.03){\textcolor{red}{\line(0,1){0.05}}}
\thicklines
\put(1.525,5.02){\textcolor{darkgreen}{\line(0,1){0.46}}}
\put(0.245,5.02){\textcolor{darkgreen}{\line(0,1){0.46}}}
\put(0.295,5.02){\textcolor{darkgreen}{\line(0,1){0.46}}}
\put(0.975,5.02){\textcolor{darkgreen}{\line(0,1){0.46}}}
\put(1.025,5.02){\textcolor{darkgreen}{\line(0,1){0.46}}}
\put(1.475,5.02){\textcolor{darkgreen}{\line(0,1){0.46}}}
\put(0.02,5.42){\textcolor{darkgreen}{\rule{0.23cm}{0.08cm}}}
\put(0.02,5.22){\textcolor{darkgreen}{\rule{0.23cm}{0.08cm}}}
\put(0.29,5.225){\textcolor{darkgreen}{\rule{0.69cm}{0.17cm}}}
\put(0.29,5.02){\textcolor{darkgreen}{\rule{0.69cm}{0.08cm}}}
\put(1.52,5.32){\textcolor{darkgreen}{\rule{0.46cm}{0.16cm}}}
\put(1.52,5.02){\textcolor{darkgreen}{\rule{0.46cm}{0.08cm}}}
\thicklines
\put(3.53,6.97){\textcolor{darkblue}{\line(1,0){1.46}}}
\put(3.53,6.53){\textcolor{darkblue}{\line(1,0){1.46}}}
\put(3.53,6.53){\textcolor{darkblue}{\line(0,1){0.46}}}
\put(4.97,6.53){\textcolor{darkblue}{\line(0,1){0.46}}}
\put(3.53,8.21){\textcolor{darkblue}{\line(1,0){1.46}}}
\put(3.53,7.53){\textcolor{darkblue}{\line(1,0){1.46}}}
\put(3.53,7.53){\textcolor{darkblue}{\line(0,1){0.69}}}
\put(4.97,7.53){\textcolor{darkblue}{\line(0,1){0.69}}}
\put(3.53,8.48){\textcolor{darkblue}{\line(1,0){1.46}}}
\put(3.53,8.26){\textcolor{darkblue}{\line(1,0){1.46}}}
\put(3.53,8.26){\textcolor{darkblue}{\line(0,1){0.23}}}
\put(4.97,8.26){\textcolor{darkblue}{\line(0,1){0.23}}}
\put(3.80,8.25){\textcolor{darkblue}{\rule{0.16cm}{0.23cm}}}
\put(4.18,8.25){\textcolor{darkblue}{\rule{0.08cm}{0.23cm}}}
\put(3.92,7.52){\textcolor{darkblue}{\rule{0.16cm}{0.69cm}}}
\put(4.62,7.52){\textcolor{darkblue}{\rule{0.08cm}{0.69cm}}}
\put(3.90,6.52){\textcolor{darkblue}{\rule{0.08cm}{0.46cm}}}
\put(4.40,6.52){\textcolor{darkblue}{\rule{0.08cm}{0.46cm}}}
\put(4.82,6.52){\textcolor{darkblue}{\rule{0.16cm}{0.46cm}}}
%
%
\put(4.10,5.63){\textcolor{red}{$A$}}
\put(3.53,5.47){\textcolor{red}{\line(1,0){1.46}}}
\put(3.53,5.03){\textcolor{red}{\line(1,0){1.46}}}
\put(3.53,5.03){\textcolor{red}{\line(0,1){0.46}}}
\put(4.97,5.03){\textcolor{red}{\line(0,1){0.46}}}
\put(3.63,5.03){\textcolor{red}{\line(0,1){0.23}}}
\put(3.93,5.13){\textcolor{red}{\line(0,1){0.12}}}
\put(4.03,5.23){\textcolor{red}{\line(0,1){0.10}}}
\put(4.43,5.23){\textcolor{red}{\line(0,1){0.20}}}
\put(4.53,5.23){\textcolor{red}{\line(0,1){0.05}}}
\put(4.73,5.03){\textcolor{red}{\line(0,1){0.05}}}
\put(4.97,5.17){$-$}
\put(5.25,5.63){\textcolor{darkgreen}{$L$}}
\put(5.22,5.48){\textcolor{darkgreen}{\line(1,0){0.23}}}
\put(5.22,5.02){\textcolor{darkgreen}{\line(1,0){0.23}}}
\put(5.22,5.02){\textcolor{darkgreen}{\line(0,1){0.46}}}
\put(5.45,5.02){\textcolor{darkgreen}{\line(0,1){0.46}}}
\put(5.22,5.22){\textcolor{darkgreen}{\rule{0.23cm}{0.08cm}}}
\put(5.22,5.42){\textcolor{darkgreen}{\rule{0.23cm}{0.08cm}}}
\put(5.97,5.63){\textcolor{darkblue}{$U$}}
\put(5.49,5.355){\textcolor{darkblue}{\line(1,0){1.46}}}
\put(5.49,5.145){\textcolor{darkblue}{\line(1,0){1.46}}}
\put(5.49,5.145){\textcolor{darkblue}{\line(0,1){0.23}}}
\put(6.95,5.145){\textcolor{darkblue}{\line(0,1){0.23}}}
\put(5.75,5.145){\textcolor{darkblue}{\rule{0.16cm}{0.23cm}}}
\put(6.15,5.145){\textcolor{darkblue}{\rule{0.08cm}{0.23cm}}}
\put(6.96,5.17){$-$}
\put(7.45,5.63){\textcolor{darkgreen}{$L$}}
\put(7.19,5.48){\textcolor{darkgreen}{\line(1,0){0.69}}}
\put(7.19,5.02){\textcolor{darkgreen}{\line(1,0){0.69}}}
\put(7.19,5.02){\textcolor{darkgreen}{\line(0,1){0.46}}}
\put(7.88,5.02){\textcolor{darkgreen}{\line(0,1){0.46}}}
\put(7.19,5.22){\textcolor{darkgreen}{\rule{0.69cm}{0.16cm}}}
\put(7.19,5.02){\textcolor{darkgreen}{\rule{0.69cm}{0.08cm}}}
\put(8.50,5.63){\textcolor{darkblue}{$U$}}
\put(7.92,5.600){\textcolor{darkblue}{\line(1,0){1.46}}}
\put(7.92,4.895){\textcolor{darkblue}{\line(1,0){1.46}}}
\put(7.92,4.895){\textcolor{darkblue}{\line(0,1){0.70}}}
\put(9.38,4.895){\textcolor{darkblue}{\line(0,1){0.70}}}
\put(8.36,4.895){\textcolor{darkblue}{\rule{0.16cm}{0.69cm}}}
\put(9.02,4.895){\textcolor{darkblue}{\rule{0.08cm}{0.69cm}}}
\put(9.37,5.17){$-$}
\put( 9.75,5.63){\textcolor{darkgreen}{$L$}}
\put( 9.62,5.48){\textcolor{darkgreen}{\line(1,0){0.46}}}
\put( 9.62,5.02){\textcolor{darkgreen}{\line(1,0){0.46}}}
\put( 9.62,5.02){\textcolor{darkgreen}{\line(0,1){0.46}}}
\put(10.08,5.02){\textcolor{darkgreen}{\line(0,1){0.46}}}
\put( 9.62,5.32){\textcolor{darkgreen}{\rule{0.46cm}{0.16cm}}}
\put( 9.62,5.02){\textcolor{darkgreen}{\rule{0.46cm}{0.08cm}}}
\put(10.60,5.63){\textcolor{darkblue}{$U$}}
\put(10.12,5.48){\textcolor{darkblue}{\line(1,0){1.50}}}
\put(10.12,5.02){\textcolor{darkblue}{\line(1,0){1.50}}}
\put(10.12,5.02){\textcolor{darkblue}{\line(0,1){0.46}}}
\put(11.58,5.02){\textcolor{darkblue}{\line(0,1){0.46}}}
\put(10.56,5.02){\textcolor{darkblue}{\rule{0.08cm}{0.46cm}}}
\put(11.08,5.02){\textcolor{darkblue}{\rule{0.08cm}{0.46cm}}}
\put(11.47,5.02){\textcolor{darkblue}{\rule{0.16cm}{0.46cm}}}
\put(3.70,4.7){\textit{copy}}
\put(4.60,4.7){$\downarrow$}
\put(3.53,4.47){\textcolor{black}{\line(1,0){1.46}}}
\put(3.53,4.03){\textcolor{black}{\line(1,0){1.46}}}
\put(3.53,4.03){\textcolor{black}{\line(0,1){0.46}}}
\put(4.97,4.03){\textcolor{black}{\line(0,1){0.46}}}
\put(3.63,4.03){\textcolor{black}{\line(0,1){0.23}}}
\put(3.93,4.13){\textcolor{black}{\line(0,1){0.12}}}
\put(4.03,4.23){\textcolor{black}{\line(0,1){0.10}}}
\put(4.43,4.23){\textcolor{black}{\line(0,1){0.20}}}
\put(4.53,4.23){\textcolor{black}{\line(0,1){0.05}}}
\put(4.73,4.03){\textcolor{black}{\line(0,1){0.05}}}
\put(3.50,3.4){\textit{gather}}
\put(4.50,3.4){$\searrow$}
\put(3.82,4.17){\textcolor{red}{\line(1,0){.16}}}\put(3.82,4.25){\textcolor{red}{\line(1,0){.16}}}\put(3.82,4.17){\textcolor{red}{\line(0,1){.08}}}\put(3.98,4.17){\textcolor{red}{\line(0,1){.08}}}
\put(4.22,4.17){\textcolor{red}{\line(1,0){.08}}}\put(4.22,4.25){\textcolor{red}{\line(1,0){.08}}}\put(4.22,4.17){\textcolor{red}{\line(0,1){.08}}}\put(4.30,4.17){\textcolor{red}{\line(0,1){.08}}}
\put(3.82,4.37){\textcolor{red}{\line(1,0){.16}}}\put(3.82,4.45){\textcolor{red}{\line(1,0){.16}}}\put(3.82,4.37){\textcolor{red}{\line(0,1){.08}}}\put(3.98,4.37){\textcolor{red}{\line(0,1){.08}}}
\put(4.22,4.37){\textcolor{red}{\line(1,0){.08}}}\put(4.22,4.45){\textcolor{red}{\line(1,0){.08}}}\put(4.22,4.37){\textcolor{red}{\line(0,1){.08}}}\put(4.30,4.37){\textcolor{red}{\line(0,1){.08}}}
%
%
%
\put(2.60,2.65){\textit{GEMM}}
\put(4.02,2.82){\textcolor{red}{\rule{0.16cm}{0.08cm}}}
\put(4.23,2.82){\textcolor{red}{\rule{0.08cm}{0.08cm}}}
\put(4.02,2.72){\textcolor{red}{\rule{0.16cm}{0.08cm}}}
\put(4.23,2.72){\textcolor{red}{\rule{0.08cm}{0.08cm}}}
\put(4.33,2.715){$-$}
\put(4.62,2.82){\textcolor{darkgreen}{\rule{0.23cm}{0.08cm}}}
\put(4.62,2.72){\textcolor{darkgreen}{\rule{0.23cm}{0.08cm}}}
\put(4.97,2.69){\textcolor{darkblue}{\rule{0.16cm}{0.23cm}}}
\put(5.16,2.69){\textcolor{darkblue}{\rule{0.08cm}{0.23cm}}}
%
%
\put(5.00,1.9){\textit{scatter}}
\put(4.50,1.9){$\searrow$}
\put(4.53,1.47){\textcolor{black}{\line(1,0){1.46}}}
\put(4.53,1.03){\textcolor{black}{\line(1,0){1.46}}}
\put(4.53,1.03){\textcolor{black}{\line(0,1){0.46}}}
\put(5.97,1.03){\textcolor{black}{\line(0,1){0.46}}}
\put(4.63,1.03){\textcolor{black}{\line(0,1){0.23}}}
\put(4.93,1.13){\textcolor{black}{\line(0,1){0.12}}}
\put(5.03,1.23){\textcolor{black}{\line(0,1){0.10}}}
\put(5.43,1.23){\textcolor{black}{\line(0,1){0.20}}}
\put(5.53,1.23){\textcolor{black}{\line(0,1){0.05}}}
\put(5.73,1.03){\textcolor{black}{\line(0,1){0.05}}}
\put(4.82,1.17){\textcolor{red}{\rule{0.16cm}{0.08cm}}}
\put(5.22,1.17){\textcolor{red}{\rule{0.08cm}{0.08cm}}}
\put(4.82,1.37){\textcolor{red}{\rule{0.16cm}{0.08cm}}}
\put(5.22,1.37){\textcolor{red}{\rule{0.08cm}{0.08cm}}}
%
%
\put(5.73,4.47){\textcolor{black}{\line(1,0){1.46}}}
\put(5.73,4.03){\textcolor{black}{\line(1,0){1.46}}}
\put(5.73,4.03){\textcolor{black}{\line(0,1){0.46}}}
\put(7.17,4.03){\textcolor{black}{\line(0,1){0.46}}}
\put(5.83,4.03){\textcolor{black}{\line(0,1){0.23}}}
\put(6.13,4.13){\textcolor{black}{\line(0,1){0.12}}}
\put(6.23,4.23){\textcolor{black}{\line(0,1){0.10}}}
\put(6.63,4.23){\textcolor{black}{\line(0,1){0.20}}}
\put(6.73,4.23){\textcolor{black}{\line(0,1){0.05}}}
\put(6.93,4.03){\textcolor{black}{\line(0,1){0.05}}}
\put(6.02,4.41){\rule{0.16cm}{0.08cm}}
\put(6.42,4.41){\rule{0.08cm}{0.08cm}}
\put(6.02,4.21){\rule{0.16cm}{0.08cm}}
\put(6.42,4.21){\rule{0.08cm}{0.08cm}}
\put(5.70,3.4){\textit{gather}}
\put(6.70,3.4){$\searrow$}
\put(6.18,4.22){\textcolor{red}{\line(1,0){.16}}}\put(6.18,4.38){\textcolor{red}{\line(1,0){.16}}}\put(6.18,4.22){\textcolor{red}{\line(0,1){.16}}}\put(6.34,4.22){\textcolor{red}{\line(0,1){.16}}}
\put(6.18,4.03){\textcolor{red}{\line(1,0){.16}}}\put(6.18,4.11){\textcolor{red}{\line(1,0){.16}}}\put(6.18,4.03){\textcolor{red}{\line(0,1){.08}}}\put(6.34,4.03){\textcolor{red}{\line(0,1){.08}}}
\put(6.78,4.22){\textcolor{red}{\line(1,0){.08}}}\put(6.78,4.38){\textcolor{red}{\line(1,0){.08}}}\put(6.78,4.22){\textcolor{red}{\line(0,1){.16}}}\put(6.86,4.22){\textcolor{red}{\line(0,1){.16}}}
\put(6.78,4.03){\textcolor{red}{\line(1,0){.08}}}\put(6.78,4.11){\textcolor{red}{\line(1,0){.08}}}\put(6.78,4.03){\textcolor{red}{\line(0,1){.08}}}\put(6.86,4.03){\textcolor{red}{\line(0,1){.08}}}
%
%
%
\put(6.22,2.74){\textcolor{red}{\rule{0.16cm}{0.16cm}}}
\put(6.43,2.74){\textcolor{red}{\rule{0.08cm}{0.16cm}}}
\put(6.22,2.62){\textcolor{red}{\rule{0.16cm}{0.08cm}}}
\put(6.43,2.62){\textcolor{red}{\rule{0.08cm}{0.08cm}}}
\put(6.53,2.70){$-$}
\put(6.82,2.74){\textcolor{darkgreen}{\rule{0.69cm}{0.16cm}}}
\put(6.82,2.62){\textcolor{darkgreen}{\rule{0.69cm}{0.08cm}}}
\put(7.61,2.40){\textcolor{darkblue}{\rule{0.16cm}{0.69cm}}}
\put(7.82,2.40){\textcolor{darkblue}{\rule{0.08cm}{0.69cm}}}
%
%
\put(7.20,1.9){\textit{scatter}}
\put(6.70,1.9){$\searrow$}
\put(6.73,1.47){\textcolor{black}{\line(1,0){1.46}}}
\put(6.73,1.03){\textcolor{black}{\line(1,0){1.46}}}
\put(6.73,1.03){\textcolor{black}{\line(0,1){0.46}}}
\put(8.17,1.03){\textcolor{black}{\line(0,1){0.46}}}
\put(6.83,1.03){\textcolor{black}{\line(0,1){0.23}}}
\put(7.13,1.13){\textcolor{black}{\line(0,1){0.12}}}
\put(7.23,1.23){\textcolor{black}{\line(0,1){0.10}}}
\put(7.63,1.23){\textcolor{black}{\line(0,1){0.20}}}
\put(7.73,1.23){\textcolor{black}{\line(0,1){0.05}}}
\put(7.93,1.03){\textcolor{black}{\line(0,1){0.05}}}
\put(7.02,1.22){\rule{0.16cm}{0.08cm}}
\put(7.42,1.22){\rule{0.08cm}{0.08cm}}
\put(7.02,1.40){\rule{0.16cm}{0.08cm}}
\put(7.42,1.40){\rule{0.08cm}{0.08cm}}
\put(7.18,1.22){\textcolor{red}{\rule{0.16cm}{0.16cm}}}
\put(7.18,1.03){\textcolor{red}{\rule{0.16cm}{0.08cm}}}
\put(7.78,1.22){\textcolor{red}{\rule{0.08cm}{0.16cm}}}
\put(7.78,1.03){\textcolor{red}{\rule{0.08cm}{0.08cm}}}
%
%
\put(7.93,4.47){\textcolor{black}{\line(1,0){1.46}}}
\put(7.93,4.03){\textcolor{black}{\line(1,0){1.46}}}
\put(7.93,4.03){\textcolor{black}{\line(0,1){0.46}}}
\put(9.37,4.03){\textcolor{black}{\line(0,1){0.46}}}
\put(8.03,4.03){\textcolor{black}{\line(0,1){0.23}}}
\put(8.33,4.13){\textcolor{black}{\line(0,1){0.12}}}
\put(8.43,4.23){\textcolor{black}{\line(0,1){0.10}}}
\put(8.83,4.23){\textcolor{black}{\line(0,1){0.20}}}
\put(8.93,4.23){\textcolor{black}{\line(0,1){0.05}}}
\put(9.13,4.03){\textcolor{black}{\line(0,1){0.05}}}
\put(8.62,4.41){\rule{0.08cm}{0.08cm}}
\put(8.22,4.41){\rule{0.16cm}{0.08cm}}
\put(8.62,4.21){\rule{0.08cm}{0.08cm}}
\put(8.22,4.21){\rule{0.16cm}{0.08cm}}
\put(8.98,4.22){\rule{0.08cm}{0.16cm}}
\put(8.38,4.22){\rule{0.16cm}{0.16cm}}
\put(8.38,4.03){\rule{0.16cm}{0.08cm}}
\put(8.98,4.03){\rule{0.08cm}{0.08cm}}
\put(7.90,3.4){\textit{gather}}
\put(8.90,3.4){$\searrow$}
\put(8.38,4.30){\textcolor{red}{\line(1,0){.08}}}\put(8.38,4.46){\textcolor{red}{\line(1,0){.08}}}\put(8.38,4.30){\textcolor{red}{\line(0,1){.16}}}\put(8.46,4.30){\textcolor{red}{\line(0,1){.16}}}
\put(8.38,4.03){\textcolor{red}{\line(1,0){.08}}}\put(8.38,4.11){\textcolor{red}{\line(1,0){.08}}}\put(8.38,4.03){\textcolor{red}{\line(0,1){.08}}}\put(8.46,4.03){\textcolor{red}{\line(0,1){.08}}}
\put(8.82,4.30){\textcolor{red}{\line(1,0){.08}}}\put(8.82,4.46){\textcolor{red}{\line(1,0){.08}}}\put(8.82,4.30){\textcolor{red}{\line(0,1){.16}}}\put(8.90,4.30){\textcolor{red}{\line(0,1){.16}}}
\put(8.82,4.03){\textcolor{red}{\line(1,0){.08}}}\put(8.82,4.11){\textcolor{red}{\line(1,0){.08}}}\put(8.82,4.03){\textcolor{red}{\line(0,1){.08}}}\put(8.90,4.03){\textcolor{red}{\line(0,1){.08}}}
\put(9.22,4.30){\textcolor{red}{\line(1,0){.16}}}\put(9.22,4.46){\textcolor{red}{\line(1,0){.16}}}\put(9.22,4.30){\textcolor{red}{\line(0,1){.16}}}\put(9.38,4.30){\textcolor{red}{\line(0,1){.16}}}
\put(9.22,4.03){\textcolor{red}{\line(1,0){.16}}}\put(9.22,4.11){\textcolor{red}{\line(1,0){.16}}}\put(9.22,4.03){\textcolor{red}{\line(0,1){.08}}}\put(9.38,4.03){\textcolor{red}{\line(0,1){.08}}}
%
%
%
\put(8.29,2.76){\textcolor{red}{\rule{0.08cm}{0.16cm}}}
\put(8.29,2.63){\textcolor{red}{\rule{0.08cm}{0.08cm}}}
\put(8.42,2.76){\textcolor{red}{\rule{0.08cm}{0.16cm}}}
\put(8.42,2.63){\textcolor{red}{\rule{0.08cm}{0.08cm}}}
\put(8.55,2.76){\textcolor{red}{\rule{0.16cm}{0.16cm}}}
\put(8.55,2.63){\textcolor{red}{\rule{0.16cm}{0.08cm}}}
\put( 8.73,2.70){$-$}
\put( 9.02,2.74){\textcolor{darkgreen}{\rule{0.46cm}{0.16cm}}}
\put( 9.02,2.62){\textcolor{darkgreen}{\rule{0.46cm}{0.08cm}}}
\put( 9.58,2.52){\textcolor{darkblue}{\rule{0.08cm}{0.46cm}}}
\put( 9.71,2.52){\textcolor{darkblue}{\rule{0.08cm}{0.46cm}}}
\put( 9.84,2.52){\textcolor{darkblue}{\rule{0.16cm}{0.46cm}}}
%
%
\put(9.40,1.9){\textit{scatter}}
\put(8.90,1.9){$\searrow$}
\put(8.93,1.47){\textcolor{black}{\line(1,0){1.46}}}
\put(8.93,1.03){\textcolor{black}{\line(1,0){1.46}}}
\put(8.93,1.03){\textcolor{black}{\line(0,1){0.46}}}
\put(10.37,1.03){\textcolor{black}{\line(0,1){0.46}}}
\put(9.03,1.03){\textcolor{black}{\line(0,1){0.23}}}
\put(9.33,1.13){\textcolor{black}{\line(0,1){0.12}}}
\put(9.43,1.23){\textcolor{black}{\line(0,1){0.10}}}
\put(9.83,1.23){\textcolor{black}{\line(0,1){0.20}}}
\put(9.93,1.23){\textcolor{black}{\line(0,1){0.05}}}
\put(10.13,1.03){\textcolor{black}{\line(0,1){0.05}}}
\put(9.22,1.22){\rule{0.16cm}{0.08cm}}
\put(9.62,1.22){\rule{0.08cm}{0.08cm}}
\put(9.22,1.40){\rule{0.16cm}{0.08cm}}
\put(9.62,1.40){\rule{0.08cm}{0.08cm}}
\put(9.38,1.22){\rule{0.16cm}{0.16cm}}
\put(9.38,1.03){\rule{0.16cm}{0.08cm}}
\put(9.98,1.22){\rule{0.08cm}{0.16cm}}
\put(9.98,1.03){\rule{0.08cm}{0.08cm}}
\put(9.38,1.30){\textcolor{red}{\rule{0.08cm}{0.16cm}}}
\put(9.38,1.03){\textcolor{red}{\rule{0.08cm}{0.08cm}}}
\put(9.82,1.30){\textcolor{red}{\rule{0.08cm}{0.16cm}}}
\put(9.82,1.03){\textcolor{red}{\rule{0.08cm}{0.08cm}}}
\put(10.22,1.30){\textcolor{red}{\rule{0.16cm}{0.16cm}}}
\put(10.22,1.03){\textcolor{red}{\rule{0.16cm}{0.08cm}}}
\put(10.50,0.7){\textit{copy}}
\put(10.00,0.7){$\downarrow$}
\put(9.93,0.47){\textcolor{black}{\line(1,0){1.46}}}
\put(9.93,0.03){\textcolor{black}{\line(1,0){1.46}}}
\put(9.93,0.03){\textcolor{black}{\line(0,1){0.46}}}
\put(11.37,0.03){\textcolor{black}{\line(0,1){0.46}}}
\put(10.03,0.03){\textcolor{black}{\line(0,1){0.23}}}
\put(10.33,0.13){\textcolor{black}{\line(0,1){0.12}}}
\put(10.43,0.23){\textcolor{black}{\line(0,1){0.10}}}
\put(10.83,0.23){\textcolor{black}{\line(0,1){0.20}}}
\put(10.93,0.23){\textcolor{black}{\line(0,1){0.05}}}
\put(11.13,0.03){\textcolor{black}{\line(0,1){0.05}}}
\put(10.22,0.40){\rule{0.16cm}{0.08cm}}
\put(10.22,0.22){\rule{0.16cm}{0.08cm}}
\put(10.62,0.40){\rule{0.08cm}{0.08cm}}
\put(10.62,0.22){\rule{0.08cm}{0.08cm}}
\put(10.38,0.22){\rule{0.16cm}{0.16cm}}
\put(10.38,0.03){\rule{0.16cm}{0.08cm}}
\put(10.98,0.22){\rule{0.08cm}{0.16cm}}
\put(10.98,0.03){\rule{0.08cm}{0.08cm}}
\put(10.38,0.30){\textcolor{black}{\rule{0.08cm}{0.16cm}}}
\put(10.38,0.03){\textcolor{black}{\rule{0.08cm}{0.08cm}}}
\put(10.82,0.30){\textcolor{black}{\rule{0.08cm}{0.16cm}}}
\put(10.82,0.03){\textcolor{black}{\rule{0.08cm}{0.08cm}}}
\put(11.22,0.30){\textcolor{black}{\rule{0.16cm}{0.16cm}}}
\put(11.22,0.03){\textcolor{black}{\rule{0.16cm}{0.08cm}}}
\put(11.50,0.15){$\to U$}
\end{picture}
}
\end{center}
\caption{Sketch of a \rsub{level--3}{level-3} BLAS update of a block row of $U$.}
\label{ldu-l3B}
\end{figure}

The hybrid structure of $L$ and $U$ allows \rdel{easily} to store
\radd{easily} the nonzeros of one block column of $L$ in \radd{a} single dense sub\rdel{--}diagonal block
and similarly, the nonzero columns of a block row of $U$ (see Figure~\ref{ldu-mem} for a sketch).
This way each block column of $L$ only consists of one dense block
and an index array referring to the nonzero row indices of $L$.
This kind of block structure is quite analogous to the structures that
are used in supernodal sparse direct $LU$ factorization methods.
Likewise, the update of a single block column of $L$ can be computed
using dense matrix kernels based on \rsub{level--3}{level-3} BLAS.
To do so, one initially has to load the associated scalar sparse columns of $A$ into
block column buffer and for each update we first need to gather the associated
submatrices required for a \rsub{level--3}{level-3} BLAS update, perform the dense 
\rsub{matrix--matrix}{matrix-matrix} multiplication (GEMM) and then to scatter the submatrix back
to the buffer. The same procedure is repeated for $U$.
The update procedure is sketched in Figure~\ref{ldu-l3B}.

In total this leads to the basic algorithm of the block incomplete $LU$
decomposition (BILU). \rdel{We make some further, minor remarks.}
\rsub{Rather than computing $A\approx LU$ such that $L$ is unit
block lower triangular and $U$ is block upper triangular, we eventually 
compute}{Technically we compute}
$A\approx LD^{-1}U$, where $L$ and $U^T$ are unit block lower triangular
and $D$ is block diagonal.
\rsub{
On one hand this requires not only
factorizing the diagonal block during the incomplete block
factorization using dense matrix kernels (LAPACK), but it further
requires us to invert the associated diagonal blocks, again using dense
matrix kernels. This simplifies dropping as well as the 
forward/backward solve during the triangular solve inside a Krylov subspace method
since it benefits from multiplying with $D$ instead of solving a sequence
of dense linear systems in each step.
Besides, dropping $l_{ik}$ whenever 
$\|l_{ik}l_{kk}^{-1}\|$ is simplified since the inverses $l_{kk}$ are precisely
these diagonal blocks of $D$ and therefore these inverses are explicitly
available.
}{
This requires factorizing and inverting the diagonal blocks using dense matrix
kernels (LAPACK) but simplifies dropping as well as the forward/backward substitution
inside a Krylov subspace method.
}
We finally note that like the scalar \rsub{Crout--type}{Crout-type} ILU our approach
does not incorporate pivoting except inside the diagonal blocks where dense
matrix kernels \radd{based on LAPACK} are used. This is certainly a drawback, however\radd{,} as
we will demonstrate
in the section on numerical results, using a combination of several
approaches (in particular maximum weight matching, blocking strategies) 
we are still able to efficiently solve a large number of systems arising from 
practical application problems.


\section{Setting Up and Improving the Block Structures}\label{sect:preprocessing}
We will now discuss several strategies that are essential to make the
\rsub{block}{B}ILU approach from the previous section efficient. We start with some
well--established scaling and permutation strategy to improve the block
diagonal dominance. Then we use an algorithm to detect 
block structures
of the initial matrix in order to group the associated rows and column together.
Based on this block partitioning \rsub{where the blocks may already have different size} we will reorder the system in order
to reduce the fill--in. In order to detect further dense blocks that will potentially
be generated by the \rsub{block}{B}ILU, 
we will perform a simplified local ILU analysis to enlarge the block size\radd{s}
of the initial block partitioning. Finally, during
the computation of the \rsub{block}{B}ILU, we allow
\rsub{to increase the block size even further}{an even further increase in the block sizes}
whenever the additional amount of fill is moderate. 
The complete portfolio of blocking strategies is inspired by the philosophy
that creating greater but fewer dense blocks is advantageous in combination
of applying dense matrix kernels\radd{,} such as \rsub{level--3}{level-3} BLAS and LAPACK\rsub{.}{, since
these are known to better exploit the cache properties
  of the underlying hardware \cite{1191}; however, to avoid higher computational complexity,
  the maximum block size should be limited.}

\subsection{Maximum Weight Matching}\label{sect:mc64}
Unless our given matrix is symmetric and positive definite,
in the general (non\rdel{-})symmetric case we may 
encounter several (block) diagonal pivots of small magnitude (or even zero). 
A \rsub{well--established}{well-established} technique that often bypasses
this problem is the use of maximum weight matchings~\cite{olschowka:1996}
as \radd{an} alternative to pivoting.
\radd{The original idea is to find a {maximum weighted matching} of
  the associated {bipartite graph} where rows and columns of the matrix refer to the nodes and
  the matrix entries serve as edge weights~\cite{DufER17}. The matching is obtained by computing a maximum product
  transversal which is equivalent to maximizing the product of the absolute values
of the diagonal entries. Finding a maximum product transversal is a well--known 
linear assignment problem in operation research and combinatorial
optimization. Essentially, one has to take the negative logarithm of the entries
and minimize the sum of the potential diagonal entries.}
For large sparse systems, as discussed here,
an efficient algorithm
was first presented in~\cite{duko:99a}. 
\radd{The problem is solved by a sparse variant of the Kuhn--Munkres algorithm. Combinatorial algorithms such as MC64~\cite{duko:99a} are experimentally observed to be extremely fast, significantly faster than the incomplete factorization itself though theoretical bounds for computing matchings are somewhat worse~\cite{duko:99a}.
The algorithm returns
a permutation as well as two dual vectors from which one has to take the
exponential in order to get the desired diagonal scalings for the original matrix. }
\rdel{Its beneficial effect in combination 
with preconditioning methods has also been established
in~[Benzi et al. 2000]. Furthermore, these kind of maximum weight matchings
are also widely used in sparse direct solvers 
(cf., e.g.,~[Schenk et al. 2004]). }
For the present manuscript we will
use \rdel{maximum weight matching}{the associated permutation and the related diagonal scalings}
as \radd{the} default initial step replacing
$A$ by
\begin{equation} \label{eqn:mwm}
\hat A=D_lAD_r\Pi,
\end{equation}
where $D_l,D_r$ are real diagonal matrices and $\Pi$ is a permutation
matrix such that the entries $\hat a_{ij}$ of $\hat A$
satisfy $|\hat a_{ij}|\leqslant 1$ and $|\hat a_{ii}|=1$. 
\rdel{Combinatorial algorithms such as MC64~[Duff and Koster 1999] are experimentally observed to be extremely fast, significantly faster than the incomplete factorization itself though theoretical bounds for computing matchings are somewhat worse~[Duff and Koster 1999]. The original idea on which these nonsymmetric permutations and scalings are based is to find a maximum weighted matching of a bipartite graph maximizing the product of the absolute values of the diagonal entries. Finding a maximum weighted matching is a well-known  linear assignment problem in operation research and combinatorial analysis. Essentially, one has to take the negative logarithm of the entries and minimize the sum of the potential diagonal entries. The minimization problem is known as a linear-sum assignment problem or bipartite weighted matching problem in combinatorial optimization. The problem is solved by a sparse variant of the Kuhn--Munkres algorithm. 
The algorithm returns
a permutation as well as two dual vectors from which one has to take the
exponential in order to get the desired diagonal scalings for the original matrix. }
\radd{\cite{olschowka:1996} introduced these scalings and
permutation for reducing pivoting in Gaussian elimination of full
matrices.}
\rdel{The permuted and rescaled system can have better numerical properties; see,
e.g.,~[Benzi et al. 2000; Duff and Koster 1999; Schenk et al. 2004].}
\rdel{Olschowka and Neumaier~[1996] introduced these scalings and
permutation for reducing pivoting in Gaussian elimination of full
matrices.} \rdel{The first implementation for sparse matrix problems was
introduced by~\cite{duko:99a}.} 
\radd{Its beneficial effect in combination 
with preconditioning methods has been established
in~\cite{duko:99a,benzi:2000:phi}. Furthermore, \radd{these kind of} maximum weight matchings
are also widely used in sparse direct solvers 
(cf.\radd{,} e.g.\radd{,}~\cite{oschenkrogu:04}).}

\begin{exm}\label{exm:venkat50-1}
Throughout this paper we will use the following matrix $A$ as a guided example
to illustrate the components of our numerical method.
The matrix \textbf{venkat50} has been taken from the {SuiteSparse Matrix Collection}\footnote{\url{https://sparse.tamu.edu/}}.
Its size is $n=62424$ with $1717777$ nonzero entries. This means that
on the average the matrix has about $27.5$ nonzero entries per row/column.
The system is nonsingular and arises from the discretization of the
\radd{2-dimensional (}2D\radd{)} Euler equations and \radd{from} there \rdel{from} a solver over several time steps 
(this matrix refers to time step $50$). The performance of sparse direct solvers
is well documented in the collection. At this moment we will use this
matrix to illustrate the effect of using maximum weight matchings.
To do so, we compute for all columns their maximum in absolute values
and call it $c_{j}$, $j=1,\dots,n$\radd{,} and similarly we proceed for the
rows to obtain $r_{i}$, $i=1,\dots,n$. In Figure~\ref{venkat50-matching}
we sketch the pattern of the entries of $A$ satisfying $|a_{ij}|\geqslant0.95\cdot\min\{r_i,c_j\}$. Analogously we sketch the pattern for $\hat A=D_lAD_r\Pi$\radd{, defined according to~\eqref{eqn:mwm}}.

\begin{figure}
 \begin{center}
\includegraphics[width=0.43\textwidth]{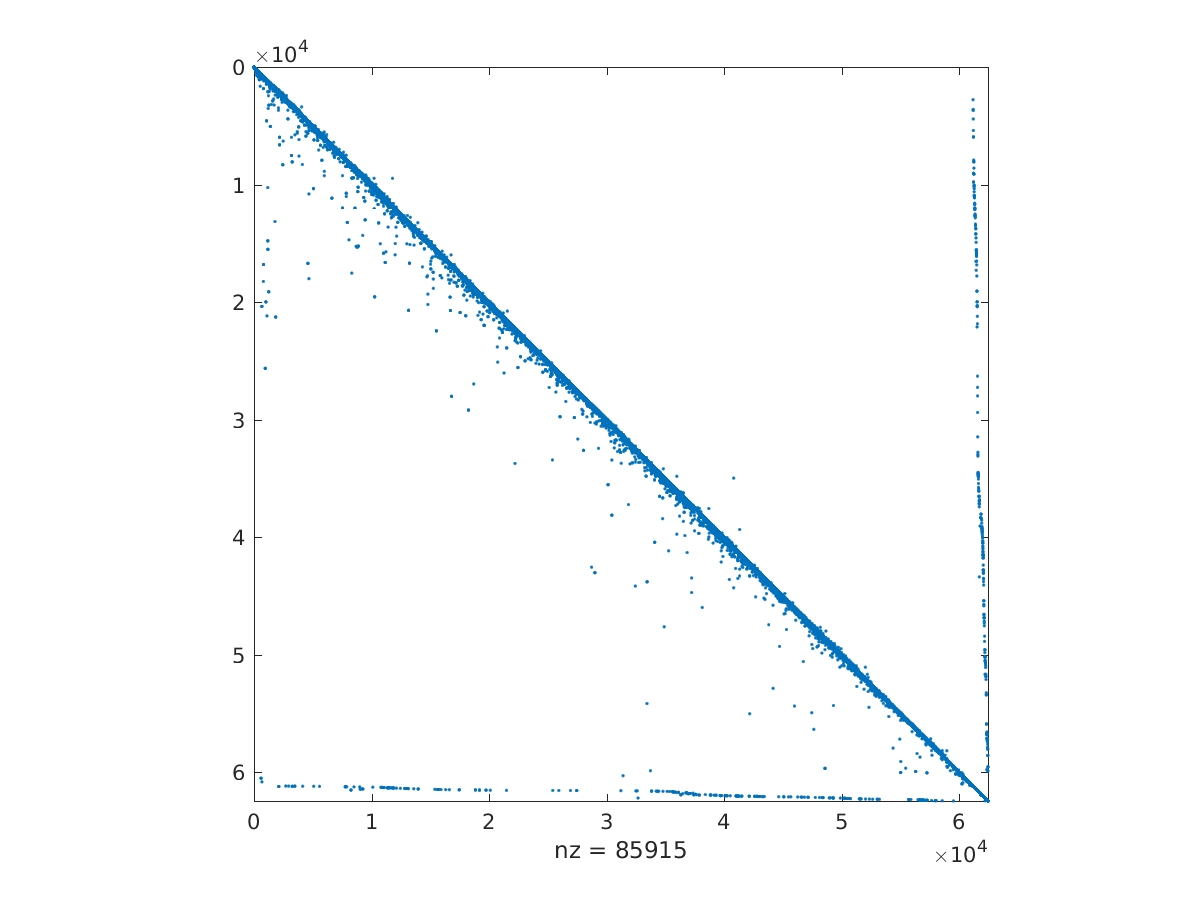} 
\qquad
\includegraphics[width=0.43\textwidth]{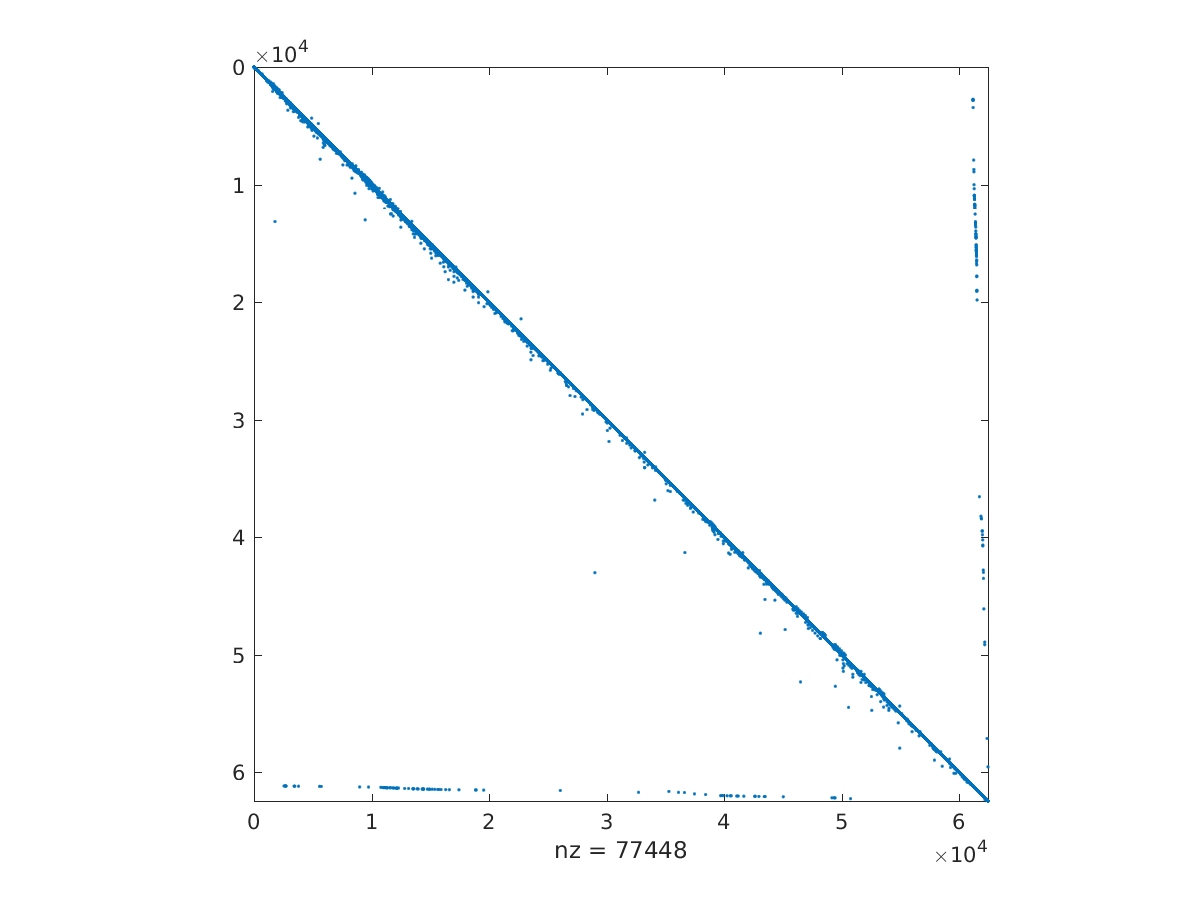} 
 \end{center}
\caption{Large entries in \rsub{modulus}{absolute value} for $A$ (left) and $\hat A$ (right).}
\label{venkat50-matching}
\end{figure}
As we can observe from Figure~\ref{venkat50-matching}, the preprocessed matrix
has significantly \rsub{less}{fewer} large entries than the original matrix and, most of them
are \rdel{not only} on the main diagonal \rsub{but also}{or at least} close to it.
\end{exm}

Since the \rdel{block incomplete factorization (}BILU\rsub{)}{,} as well
as its scalar counter part as far as discussed in this paper,
do not use further pivoting\rsub{,}{;} the use of maximum weight matching
is an essential generic preprocessing step 
in obtaining a successful incomplete factorization in
many test cases though there are certainly some practical problems where 
the use of maximum weight matchings is less beneficial.


\subsection{Cosine-Based Preprocessing}\label{sect:cosine}
The initial preprocessing step using maximum weight
matching\radd{, it is hoped, } simply improves the scalar diagonal dominance\rdel{, hopefully}.
\rsub{So far we have not discussed a strategy to set up an initial block
structure of the matrix that can possibly be improved during the
approximate factorization.
Here we propose as initial option to use the cosine--based
strategy as suggested in~[Saad 2003a]. I}{We now propose a cosine-based strategy to initialize a block structure of the matrix and that could be possibly improved during the approximate factorization process~\cite{Saa03a}; i}n our numerical
experiments we will use BILU with and without the cosine-based
blocking to illustrate the overall performance of the code.
\rdel{We will briefly explain the main idea of the cosine-based blocking
strategy.}
Given two rows $a_i^T=e_i^T\hat A$, $a_j^T=e_j^T\hat A$ 
of a matrix $\hat A$, their nonzero pattern can be represented by two
row vectors $c_i^T$, $c_j^T$ which have values $1$ if and only if the associated
entries of $a_i$, $a_j$ are nonzero and $0$ otherwise.
The major observation is that two rows 
have almost the same nonzero pattern if their formal scalar product satisfies
$c_i^Tc_j\approx \|c_i\|\cdot \|c_j\|$. Since this computation is integer\rsub{--}{ }based,
a simple counting strategy for $nz(a_i \cap a_j)^2\geqslant \tau\cdot nz(a_i)\cdot nz(a_j)$ is \rsub{enough}{sufficient}, where $\tau\in[0,1]$ is a prescribed threshold.
In~\cite{Saa03a}, $\tau=0.8$ is suggested which we will use as well.
The algorithm uses the (sparse) 
pattern of the upper triangular part of $\hat A\hat A^T$ as long
as the associated indices are not yet associated with some \radd{diagonal} block.

Overall, whenever we use the \rsub{cosine--based}{cosine-based} strategy, we replace
$\hat A$ by
\begin{equation} \label{eqn:cos}
\tilde A = Q^T\hat A Q,
\end{equation} 
where $Q$ is \rsub{a permutation matrix grouping together columns and rows
such that $\tilde A$ has an improve block pattern}{the permutation matrix generated by the cosine-based approach grouping together columns and rows of $\hat{A}$. This results in $\tilde{A}$ having an improved block pattern}.
We finally like to mention that beside its benefits, the cosine
strategy may become extremely inefficient \rsub{when}{for cases where} $\hat A\hat A^T$ becomes relative\radd{ly}
dense although $\hat A$ is relatively sparse. \rsub{This could be, e.g., the
case when $\hat A$ has a few columns/rows with relatively many
nonzero entries.}{This situation might verify, e.g., when some of the rows of $\hat{A}$, even though in limited number, are densely populated with nonzeros.}
For this reason we use a slightly modified version of Saad's cosine
strategy that ignores rows/columns having too many nonzero entries.
This is done by a simple statistical argument computing the average 
number $\mu$ of nonzeros per row/column as well as the associated
standard deviations $\sigma_{r,c}$. Rows (resp.\radd{,} columns) exceeding
$\mu+2\sigma_{r,c}$ nonzeros are ignored for the cosine blocking strategy.

\begin{exm}\label{exm:venkat50-2}
We continue Example~\ref{exm:venkat50-1} and illustrate\radd{,} for the matrix
$\hat A$ obtained after maximum weight matching has been applied, how many
blocks
were detected by the \rsub{cosine--based}{cosine-based} method. Here our results already refer to the
modified version\rsub{.}{:}

\begin{center}
\begin{tabular}{ccccccccc}
system size && \# dgl. blocks && max. size && avg. size && std. deviation\\
\cline{1-1}\cline{3-3}\cline{5-5}\cline{7-7}\cline{9-9}
$62424$    && $15723$          && $4$                && $3.97$  && $0.297$
\end{tabular}
\end{center}

We can see that the majority of blocks detected by the cosine algorithm ha\rsub{s}{ve} a block
size $4$.
%
\end{exm}


\subsection{Symmetric Reordering}\label{sect:metis}
After having identified potential initial blocks using the
cosine-based strategy (or even when leaving it out), we next will 
reorder the system $\tilde A$ respecting the given block pattern.
If the cosine strategy was not used, we \radd{would} simply use the scalar
partitioning instead\radd{, i.e., the original matrix}.
However, replacing $\tilde A$ by its companion
matrix that compresses each block of $\tilde A$ into a scalar in a
straightforward manner, we may reorder the associated companion matrix
$B$ using standard symmetric reordering strategies such as approximate
minimum degree~\cite{AmeDD96} or nested dissection~\cite{karypis:98,LasK13}
as implemented in the METIS package.
Here\radd{,} for simplicity\radd{,} we restrict ourselves to the nested dissection ordering
as implemented in METIS in order to reduce the fill--in further.
\rsub{Whenever}{After} METIS is applied to the compressed companion matrix $B$
we expand the associated permutation matrix $P_B$ \rsub{again in order}{to a block permutation matrix $P$ that preserves the block structure of $A$ and thus} obtain a reordering for the original matrix $\tilde A$ respecting the
block partitioning. This gives a symmetrically reordered matrix
$\check{A}$, where
\begin{equation} \label{eqn:ordering}
\check{A} = P^T \tilde A P.
\end{equation}
We sketch the approach in the following illustration\rsub{.}{:}
\[
\underbrace{\tiny\left(
\begin{array}{@{}c@{\,}c|c@{\,}c|c@{\,}c@{}}
\ast &\ast &\ast &\ast &\\
\ast &\ast &     & \ast &\\\hline
\ast &\ast &\ast &      &\ast &\ast \\
     &\ast &\ast &\ast &\ast &\\\hline
     &     &\ast &     &\ast &\ast \\
     &     &\ast &\ast &     &\ast
\end{array}
\right)}_{\hat A}
 \stackrel{compress}{\to} 
\underbrace{\tiny\left(
\begin{array}{@{}c@{\,}c@{\,}c@{}}
\ast & \ast & \\
\ast & \ast & \ast \\
     & \ast & \ast  
 \end{array}
\right)}_B \stackrel{reorder}{\to} 
\underbrace{\tiny\left(
\begin{array}{@{}c@{\,}c@{\,}c@{}}
\ast  &     & \ast\\
      & \ast& \ast\\  
\ast  & \ast& \ast 
 \end{array}
\right)}_{P_B^TBP_B} \stackrel{expand}{\to} 
\underbrace{\tiny\left(
\begin{array}{@{}c@{\,}c|c@{\,}c|c@{\,}c@{}}
\ast &\ast &     &    &\ast &\ast  \\
\ast &\ast &     &    &     & \ast \\\hline
     &     &\ast &\ast&\ast &      \\
     &     &     &\ast&\ast &\ast  \\\hline
\ast &\ast &\ast &\ast&\ast &      \\
     &\ast &\ast &    &\ast &\ast 
\end{array}
\right)}_{P^T\hat AP}.
\]

\begin{exm}\label{exm:venkat50-3}
Finally we illustrate in Figure~\ref{venkat50-matching-metis}
how the matrix from Examples~\ref{exm:venkat50-1}
and~\ref{exm:venkat50-2} is reordered with nested dissection following
an initial blocking strategy obtained by the cosine algorithm.

\begin{figure}
 \begin{center}
\includegraphics[width=0.43\textwidth]{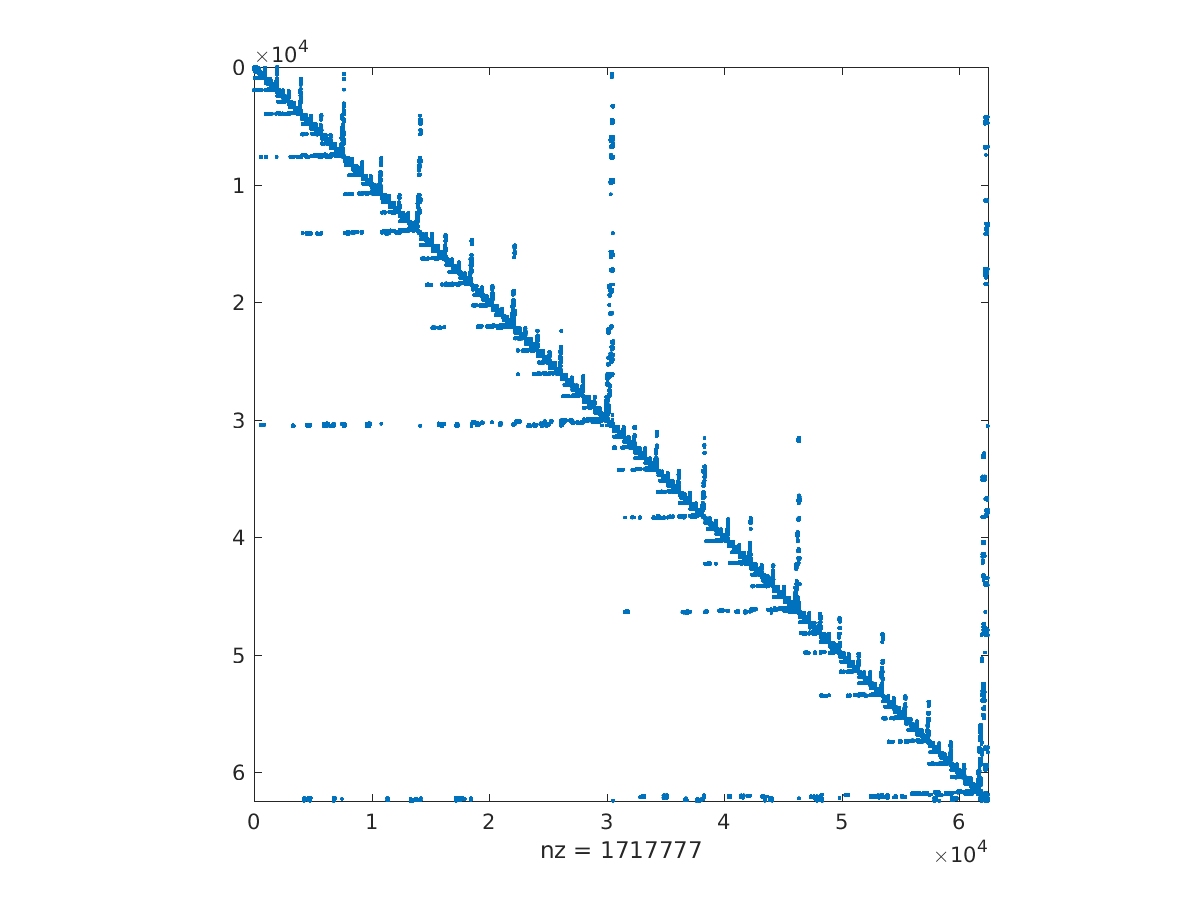} 
 \end{center}
\caption{Matrix \textbf{venkat50} after maximum weight matching, \rsub{cosine--based}{cosine-based} blocking and METIS reordering based on the compressed graph}
\label{venkat50-matching-metis}
\end{figure}

Given this preprocessed matrix we could now start with the \rsub{Crout--type block
}{Crout-type B}ILU simply inheriting the block structure and its reordering. \rsub{Beside}{Apart from computing} the 
incomplete factorization we need to solve the associated linear system
iteratively. Here we use\radd{,} for simplicity\radd{,} the restarted GMRES~\cite{SaaS86} method with
restart length $30$ until the relative residual is reduced by $10^{-6}$.
As right\radd{-}hand side we use the vector $b$ with all ones. 
In our experimental environment, the code was implemented in C but using
a CMEX interface to MATLAB (R2015b). 
The same applies to the forward/backward
solve. 
This MATLAB release uses 
Intel MKL 11.1.1 including BLAS and LAPACK 3.4.1.
The results were \rsub{carried out}{obtained} on a 
single node with 1 TB main memory and $4$ Intel Xeon E7-4880 v2 @ 2.5 GHz processors
each of them having $15$ cores on a socket leading to $60$ cores in total.
For this specific example, we compare \rdel{the block ILU method (``}BILU\rdel{'')} as
described so far
with \rsub{MATLAB's}{the MATLAB} \texttt{ilu} function and its option \texttt{\rdel{'}crout\rdel{'}} (referred hereafter as \rdel{``}ILUC\rdel{''}) which perfectly
fits with the block ILU as scalar counter part. Both methods
use maximum weight matching, the METIS reordering\radd{,} and a drop tolerance $\tau=10^{-2}$\rsub{.}{:}

\begin{center}
\begin{tabular}{ccccccccc}
     && time ILU[sec] && $\frac{nz(L+U)}{nz(A)}$ && time GMRES[sec] && \# steps\\
\cline{3-3}
\cline{5-5}
\cline{7-7}
\cline{9-9}
BILU && \phantom{0}1.9        && 4.3  && 3.0 && 29\\
ILUC && 20.0       && 2.7  && 4.7 && 95
\end{tabular}
\end{center}
Apparently the blocking strategy in combination with the \rsub{block }{B}ILU
outperforms its scalar \rsub{counter part}{counterpart} already by one order of magnitude.
\end{exm}
%
%
%
%


\subsection{Guessing Initial Block Pattern via a Simplified ILU}\label{sect:ilu1t}
\rsub{Now h}{H}aving improved the diagonal dominance and possibly having identified
initial dense blocks and reordering the associated companion matrix
using a fill-reducing method, we now could start
factorizing the matrix $\check A$ (cf. Example~\ref{exm:venkat50-3}). 
For reasons of efficiency it may
pay \rdel{off} to take a closer look at the given matrix $\check A$ before
starting the block incomplete factorization, in particular\radd{,} taking a look
at the block partitioning. \rsub{As the cosine--based algorithm only
identifies potential dense blocks in the initial matrix, direct methods
using the elimination tree explore the matrix pattern in order
to identify dense blocks that show up during the $LU$ factorization.}
{We underline that the dense blocks identified by the cosine-based analysis of the matrix pattern are unrelated to the ones defined by the elimination tree explore performed by the direct $LU$ factorization methods.}
Here\radd{,} within the context of an incomplete factorization, we cannot
expect symbolic strategies to give a sufficiently accurate prediction
about dense blocks, since dropping entries of small size
during the factorization will destroy most of the \rsub{graph--based}{graph-based}
information. \rsub{Instead, we propose}{We propose instead} to use a simplified incomplete
$LU$ factorization model that is drastically cheaper than the
\rsub{block }{B}ILU itself but might serve as a simple first\rsub{--}{ }order guess for 
the dense blocks that may show up during the factorization. \radd{There are several possibilities
  to compute block patterns for incomplete factorizations, e.g., exploiting masks like the level of fill or thresholds
  or both of them   \cite{Wat81,DazFT92,HysP01,HenRR08,ScoT11}.}

To be precise let us briefly recall the \rsub{level--of--fill}{level-of-fill} approach $ILU(p)$ 
(cf.~\cite{DazFT92,HysP01,HysP02,Saa03}).
Initially, we define a level function $lev_{ij}$ via  
\[
lev_{ij}:= 
\left\{
	\begin{aligned}
		0\rdel{,} 		\quad  & \radd{\text{if}\;} a_{ij}\not=0 \radd{,} \\
		\infty\rdel{,}  \quad  &\mbox{otherwise} \radd{.}
	\end{aligned}
\right.
\]
During the approximate factorization we are faced at step $k$ 
with updates of the form
\[
a_{ij}\leftarrow a_{ij}-\frac{a_{ik}\,a_{kj}}{a_{kk}}
\]
which modifies the level function to become 
\[
lev_{ij}=\min\{lev_{ij},lev_{ik}+lev_{kj}+1\}.
\]
Now the \rsub{level--of--fill}{level-of-fill} $ILU(p)$ only allows these kind of updates
whenever $lev_{ij}\leqslant p$. Otherwise the update is omitted.
If $a_{ij}\not=0$ before the update, then $lev_{ij}<\infty$ does not increase
anymore and one may update, whereas for $a_{ij}=0$, only updating
is permitted as long as $lev_{ik}+lev_{kj}+1\leqslant p$. This limits
the number of fill--entries. \rsub{The most common cases are $p=0$ and $p=1$.}{Often enough, smaller
  numbers of $p$ are used.}
For $p=0$, the \rsub{level--of--fill}{level-of-fill} ILU simply inherits the original pattern of
$A$ and disregards any fill outside the initial pattern.
For $p=1$, additional nonzero entries $a_{ij}$ can only be created by an
update where we have $a_{ik}\not=0$ and $a_{kj}\not=0$ in the initial
matrix $A$. In short\rsub{:}{,} \rsub{fill--in}{fill-in} is permitted by original entries only, 
but not by \rsub{fill--in}{fill-in} entries.

Another ILU approach consists of dropping entries of small size $\tau$
(referred to as $ILU(\tau)$), i.e.,
at step $k$ of the ILU
we discard $a_{ik}$, $a_{kj}$, whenever $|a_{ik}|\leqslant \tau |a_{kk}|$
(resp.\radd{,} $|a_{kj}|\leqslant \tau |a_{kk}|$). This is applied \rsub{no matter}{regardless of} whether
the entries were originally nonzero or created as \rsub{fill--in}{fill-in}.
Since we are using the Crout-based ILU (see Algorithm~\ref{iluc}),
at step $k$ only column $k$ and row $k$ of the incomplete factorization
are computed, i.e., 
$a_{ik},i\geqslant k$\radd{,} and $a_{kj},j\geqslant k$. This is why
we relate their size with respect to $|a_{kk}|$.

Certainly one could easily combine $ILU(p)$ and $ILU(\tau)$ to obtain
some kind of $ILU(p,\tau)$. For simulation \rdel{of} particularly \radd{of} $ILU(1,\tau)$,
we apply the method only locally estimating the fill pattern at step $k$.
The idea is \radd{to} simulate the behavio\rdel{u}r of the Crout ILU quickly and find
from this quick and simple simulation a good initial guess for the block
pattern.

Suppose that we plan to estimate the fill pattern of column $k$ of $L$
and row $k$ of $U$. The initial index sets 
$I=\{i|\, i\geqslant k, a_{ik}\not=0\}$, 
$J=\{j|\, j\geqslant k, a_{kj}\not=0\}$ 
consist of the associated pattern of $A$. We do not plan to update the
diagonal part $a_{kk}$ as part of the incomplete factorization process and
\radd{will} simply use $\tau\cdot |a_{kk}|$ as \radd{the} approximation, i.e,
$I$ is reduced to those indices $i$ satisfying $|a_{ik}|\geqslant \tau\cdot |a_{kk}|$. We call this set $\hat I$ (resp.\radd{,} $\hat J$ for the upper triangular part).

Next, in accordance \rsub{to}{with} the $ILU(1)$ philosophy to only allow fill-in from
the original nonzero entries, we are seeking for all nonzero entries
$a_{ij}$, $a_{jk}$  of $A$ such that $j<k<i$. These can be easily
obtained \radd{by} first looking at the nonzero pattern of column $k$ of $A$ and,
inside column $k$ only for those $j$ satisfying $j<k$. Let us denote
this set by $J_k$.
For all these indices $j\in J_k$ we need to check column $j$ of $A$
for indices $i>k$. Given $i,j$, we can easily simulate the \rsub{fill--in}{fill-in} situation
only adding $i$ to $\hat I$, if $|a_{ij}a_{jk}|\geqslant \tau |a_{jj}a_{kk}|$
is fulfilled. This refers to a \rsub{fill--in}{fill-in} $a_{ik}=-\frac{a_{ij}a_{jk}}{a_{jj}}$
which is large enough compared with $|a_{kk}|$. Again, 
$a_{ij}$, $a_{jk}$\radd{,} and $a_{jj}$ from
the original matrix $A$ are used to substitute the values of the incomplete
factorization. Thus this \radd{is only a local} analysis\rdel{ is local analysis only}.
\rsub{T}{We denote by $I$ t}he column index set we obtain \radd{by} including this type of \rsub{fill--in we will denote
by $\tilde I$}{fill-in}. We proceed similarly to obtain $\tilde J$.
As we have now computed the estimated patterns of column $k$ of
 $L$ and row $k$ of $U$, we could continue to compute similar patterns
in steps $k+1,k+2,k+3,\dots$ leading to a sequence of patterns
$\tilde I_k,\tilde I_{k+1},\tilde I_{k+2},\dots$ and
$\tilde J_k,\tilde J_{k+1},\tilde J_{k+2},\dots$.
For aggregating scalar columns/rows to build blocks we simply need to
build their union 
$\check I=\tilde I_k\cup\tilde I_{k+1}\cup\tilde I_{k+2}\cdots$ and
$\check J=\tilde J_k\cup\tilde J_{k+1}\cup\tilde J_{k+2}\cdots$ 
measuring the additional zero entries when incorporating the next
column/row and removing the entries that refer to the diagonal block,
which are considered to be part of a dense diagonal block.
This way, we can exactly compute the additional zero entries to fill
up the blocks when adding a new column/row into the current block
(certainly assuming that our local ILU analysis is accurate enough).
Suppose that this way we proceed from step $k$ to $k+l-1$ and assume
that the subdiagonal block of $L$ consists of $r$ nonzero rows 
and the superdiagonal block of $U$ consists of $s$ nonzero columns
whereas the additional zero entries are given by some number
$c$. 
This way we have $f_l=(r+s+l)\cdot l$ nonzeros in the block case
whereas the scalar case would only have $f_l-c$
nonzero entries. Going from step $k+l-1$ to step $k+l$ we obtain new
values $r',s'$ for the off--diagonal blocks and 
$f_{l+1}=(r'+s'+l+1)\cdot (l+1)$. The new scalar fill would become
$c'=f_l-c+|I_{k+l}|+|J_{k+l}|+1$.
In order to avoid an inflation of
additional wasted zero entries
 we allow \radd{one} to incorporate the $(l+1)$\rsub{.}{-st} column/row also,
as long as
\[
f_{l+1}\leqslant \frac43 \cdot c' \mbox{ or }
f_{l+1}\leqslant c' + 4 \cdot (l+1)
\]
holds.
This allows to increase the overhead of wasted zeros by $1/3$
with respect to the scalar situation or alternatively 
to have, say, $2$ additional rows in $L$ and $2$ additional columns in $U$ 
in the block partitioning 
(or 4 additional rows in $L$ but no additional column in $U$, etc\radd{.}).

In our test cases in \rsub{S}{s}ection~\ref{sect:exp}
we will illustrate how the algorithms perform
with and without the $ILU(1,\tau)$ strategy. Beside\radd{s,} we will also
demonstrate the behavio\rdel{u}r for a specific example as follows.

\begin{exm}\label{exm:venkat50-4}
We continue Examples~\ref{exm:venkat50-1}--\ref{exm:venkat50-3} and in
a first step we compare the two blocking strategies when being applied
separately and together (with the METIS reordering in\radd{-}between).
\begin{center}
\begin{tabular}{c@{}cccccccccc}
      && \# blocks && max. size && avg. size && std. deviation\\
\cline{3-3}\cline{5-5}\cline{7-7}\cline{9-9}
only cosine         && $15723$ && $4$          && $3.97$  && $0.297$\\
only $ILU(1,10^{-2})$  && $24138$ && $16$         && $2.59$  && $2.15$ \\
cosine+$ILU(1,10^{-2})$&& $13786$ && $16$         && $4.53$  && $1.50$
\end{tabular}
\end{center}
We can see that the major blocking was already obtained from the initial
cosine strategy while the $ILU(1,10^{-2})$ has added some larger blocks.
It looks as if the combination of both yields the best blocking.
Looking at the performance of the associated \rsub{block }{B}ILU variants (using the
abbreviation ``c-'' for the pure cosine strategy, ``-i'' for only using $ILU(1,10^{-2})$\radd{,}
and ``ci'' for both) we observe that  the combination of both 
blocking strategies in this example is at least comparable
with the initial blocking strategy\rsub{.}{:}
\begin{center}
\begin{tabular}{ccccccccc}
     && time ILU[sec] && $\frac{nz(L+U)}{nz(A)}$ && time GMRES[sec] && \# steps\\
\cline{3-3}
\cline{5-5}
\cline{7-7}
\cline{9-9}
BILU(c-) && 1.9           && 4.3  && 3.0 && 29\\
BILU(-i) && 3.9           && 4.0  && 4.0 && 31\\
BILU(ci) && 2.6           && 4.5  && 2.5 && 26
\end{tabular}
\end{center}
\end{exm}
%
%


\subsection{Progressive Aggregation}\label{sect:pa}
So far we have simply worked with \radd{variable} block structures that were predefined
in advance, \rdel{say} either using the cosine-based method or the ILU(1,$\tau$)
strategy or even both of them.
In order to improve the blocking further, we will merge blocks during
the factorization in the case that two consecutive block columns of 
$L$ and $U^T$ follow each other and the additional memory overhead is 
acceptable. Although this will increase the \rsub{fill--in}{fill-in} and although this
aggregation is not completely for free, the expectation is\rdel{,} that having
\rsub{less}{fewer} but larger blocks pays off in combination with the use of \rsub{level--3}{level-3} BLAS.
We note that the computation of a block column of $L$ (resp.\radd{,} block row of $U$)
in a \rsub{Crout--type}{Crout-type} \rsub{block }{B}ILU (see \rsub{S}{s}ection~\ref{sect:biluc}) 
requires \radd{one} to compute this block column at some step $k$ based
on several preceding block columns of $L$. If their number decreases but
their size increases\radd{,} while we compute the ILU, me may
expect that the \rsub{level--3}{level-3} BLAS computation leads to an acceleration \radd{as long as the maximum block size is limited to avoid that computational complexity starts dominating the process}.


Suppose that\radd{, after $k$ steps of progressive aggregation,} we have computed from our approximate $LU$ decomposition
the leading $k$ block columns/rows $L^{(k)}$, $U^{(k)}$
as well as the leading inverse block diagonal matrix $(D^{(k)})^{-1}$\rsub{,
where}{.}
\[
  \rdel{
L^{(k)}=
\left(
\begin{array}{cccc}
  I       &        &           &0  \\
L_{21}    & \ddots &           &  \\
\vdots    & \ddots &     I     &  \\
L_{k+1,1} & \cdots & L_{k+1,k} & I
\end{array}
\right), \; 
U^{(k)}=
\left(
\begin{array}{cccc}
  I       & U_{12} & \cdots    & U_{1,k+1} \\
          & \ddots & \ddots    & \vdots\\
          &        &     I     & U_{k,k+1} \\
  0       &        &           & I
\end{array}
\right),} 
\]
\[
\rdel{(D^{(k)})^{-1}=\left(
\begin{array}{ccc}
D_{11}^{-1}&        &    0       \\
           & \ddots &            \\
  0        &        & D_{kk}^{-1}
\end{array}
\right). }
\]
\rdel{Here each block column of $L^{(k)}$ is stored in a dense array, 
where only the nonzero rows are stored but it includes all columns belonging
to this block (as explained earlier in section~\ref{sect:biluc}).
A similar storage scheme is used for $(U^{(k)})^T$.} The number of block columns
in $L$ \radd{(resp. $L^{(k)}$)} has been predefined up to step $k$\radd{,} whereas for steps $k+1$\radd{,} and later\radd{,}
we still could change the block sizes easily, since this part has not \radd{yet}
been computed \rdel{yet}.
Merging block columns/rows $k-1$ and $k$ requires \radd{us} to rewrite the associated
matrices as
\begin{eqnarray*}
&&\left(
\begin{array}{cc}
  I         &    0      \\
L_{k,k-1}   &    I      \\
L_{k+1,k-1} & L_{k+1,k} 
\end{array}
\right)
\left(
\begin{array}{cc}
D_{k-1,k-1}^{-1}&    0       \\
  0        & D_{kk}^{-1}
\end{array}
\right)
\left(
\begin{array}{ccc}
  I       & U_{k-1,k} & U_{k-1,k+1} \\
  0       &     I     & U_{k,k+1}
\end{array}
\right)\\
&=&\left(
\begin{array}{cc}
  I         &    0      \\
  0         &    I      \\
\hat L_{k+1,k-1} & L_{k+1,k} 
\end{array}
\right)
\left(
\begin{array}{cc}
\hat D_{k-1,k-1}&\hat D_{k-1,k}\\
\hat D_{k-1,k}  &\hat D_{kk}
\end{array}
\right)^{-1}
\left(
\begin{array}{ccc}
  I       &     0 & \hat U_{k-1,k+1} \\
  0       &     I & U_{k,k+1}
\end{array}
\right).
\end{eqnarray*}
\rdel{Obviously the aggregated block column of $L^{(k)}$ is obtained
from its two previous block columns $k-1$ and $k$ by post multiplying the
block columns by}
\[
\rdel{\left(
\begin{array}{cc}
  I         &    0      \\
-L_{k,k-1}   &    I      
\end{array}
\right).}
\]
\rdel{From this it follows that 
$\hat L_{k+1,k-1} =L_{k+1,k-1} - L_{k+1,k}L_{k,k-1}$,
whereas $L_{k+1,k}$ remains unchanged. Similar arguments apply to $U^{(k)}$.}
The aggregated inverse block diagonal block of $(D^{(k)})^{-1}$ 
is adapted accordingly\radd{,}
leading to a larger dense inverse diagonal block.
Aggregating two consecutive block columns/rows typically increases
the fill\radd{-}in $\hat L_{k+1,k-1}$, $\hat U_{k-1,k+1}^T$\radd{,} and also in
$L_{k+1,k}$ and $U_{k,k+1}^T$, since the aggregated blocks
$[\hat L_{k+1,k-1}, L_{k+1,k}]$ 
need to have a common nonzero row pattern and
\rdel{and} $\left[\hat U_{k-1,k+1}\atop U_{k,k+1}\right]$ must have the same column pattern.
We allow to aggregate the blocks progressively\rdel{,} whenever the \rsub{additional
amount of memory}{memory increase} is mild. Suppose that block column $k-1$ consists of $p$ columns
and block column $k$ has $q$ columns.
The subdiagonal blocks $L_{k,k-1}$, $L_{k+1,k-1}$\radd{,} and $L_{k+1,k}$ may have $r,s,t$
nonzero rows and similarly  $U_{k-1,k}$, $U_{k-1,k+1}$\radd{,} and $U_{k,k+1}$ may have
$r',s',t'$. Then before the aggregation the number of nonzeros is given by
\[
\mu = p\,(p+r+s+r'+s')+ q\,(q+ t+t').
\]
Without explicitly computing $\hat L_{k+1,k-1}$ or $\hat U_{k-1,k+1}$,
we can easily compute the associated number of nonzero\rdel{s} rows $u$ of
$[\hat L_{k+1,k-1}, L_{k+1,k}]$ and nonzero columns $v$ of
$\left[\hat U_{k-1,k+1}\atop U_{k,k+1}\right]$ \radd{by} simply checking the
union of nonzero index sets ignoring any cancellations.
This gives 
\[
\nu=(p+q)\,(p+q+u+v)
\]
nonzero entries, where the diagonal block is always stored in dense format.
We let the algorithm aggregate block $k-1,k$ to become an enlarged
block $k-1$ whenever
$\nu\leqslant 1.2\cdot  \mu$ or $\nu\leqslant\mu + 2 (p+q)$ is satisfied.
Certainly one could vary these numbers and we do not claim that
they are ``best'' in some sense. The philosophy is to allow 20\% additional
\rsub{fill--in}{fill-in} or at least two rows/column (e.g.\radd{,} one in $L$ and one in $U$).
After checking some examples this has turned out to be an acceptable compromise
between \rsub{fill--in}{fill-in} and the size of the blocks. 

We note that the aggregation process is \rsub{checked in every}{always checked in} step $k$\radd{,} allowing
\rsub{to increase the block size}{the block sizes to increase} progressively\rsub{, i.e.,}{. Because of this,} it may happen that in \radd{all the} steps
$k,k+1,k+2,\dots,k+l$ \rdel{that} the current block is \rdel{always} aggregated with its
predecessor\radd{,} such that at step $k+l$ we only have one aggregated block\radd{, labeled $k-1$}.
Theoretically the \rsub{fill--in}{fill-in} could be drastically increased, but we did not observe
this in our practical experiments. This may be related to the fact that a fill-reducing
ordering (in our case nested dissection) was applied prior to the \rsub{block }{B}ILU computation.
Finally we note that the data structures of the \rsub{Crout--type}{Crout-type} \rsub{block }{B}ILU from 
\rsub{S}{s}ection~\ref{sect:biluc} can be adapted easily. Technically, the easiest implementation
has turned out to define block $k-1$ simply as void (block size $0$) and to
let the aggregated block become block $k$. This way, the auxiliary vectors
\texttt{L\_head}, \texttt{L\_list}, \texttt{L\_first} for $L$
and \texttt{U\_head}, \texttt{U\_list}, \texttt{U\_first} for $U$
from \rsub{S}{s}ection~\ref{sect:iluc} need not be changed at all and the void block $k-1$
quickly drops out step by step (since it is not longer needed for updates).

\begin{exm}\label{exm:venkat50-5}
We finish Examples~\ref{exm:venkat50-1}--\ref{exm:venkat50-4} 
by examining the additional benefits of the progressive aggregation.
In analogy to Example~\ref{exm:venkat50-4} we sketch the compression
rate of each single blocking strategy and their combination\rsub{.}{:}
\begin{center}
\begin{tabular}{c@{}cccccccccc}
      && \# blocks && max. size && avg. size && std. deviation\\
\cline{3-3}\cline{5-5}\cline{7-7}\cline{9-9}
cosine                && $15723$ && $4$          && $3.97$  && $0.297$\\
$ILU(1,10^{-2})$       && $24138$ && $16$         && $2.59$  && $2.15$ \\
progr. aggr.            && $32996$ && $13$         && $1.89$  && $1.59$ \\
\multicolumn{1}{p{5cm}}{cosine + $ILU(1,10^{-2})$ + progr. aggr.}&& $8259$ && $44$         && $7.56$  && $5.23$
\end{tabular}
\end{center}
As already observed earlier, the best blocking performance results from the combination
of all three methods.
Finally we compare the \rsub{block }{B}ILU method when using only one of the three blocking
strategies with the version that incorporates all strategies\rsub{.}{:}
\begin{center}
\begin{tabular}{ccccccccc}
     && time ILU[sec] && $\frac{nz(L+U)}{nz(A)}$ && time GMRES[sec] && \# steps\\
\cline{3-3}
\cline{5-5}
\cline{7-7}
\cline{9-9}
BILU(c-{}-) && 1.9           && 4.3  && 3.0 && 29\\
BILU(-i-)   && 3.9           && 4.0  && 4.0 && 31\\
BILU(-{}-p) && 3.7           && 3.5  && 7.4 && 52\\
BILU(cip)   && 1.9           && 4.9  && 2.3 && 26
\end{tabular}
\end{center}
For this example the overall performance is best using the three methods together
at the price of a slightly higher \rsub{fill--in}{fill-in}.
\end{exm}

We conclude this section noting that using progressive aggregation
\rdel{only} without an initial block strategy can become quite costly, since
the strategy may merge two consecutive block columns/rows several times
successively, increasing a scalar column/row to a block size of a few
hundred. It is clear that this can hardly be efficient\radd{,} in general\radd{,} and this is
why having some initial guess\rdel{es} for the block partitioning
prior to the progressive
aggregation is useful.


\subsection{Perturbing the Entries of the Diagonal Blocks} \label{sect:perturbation}
In the symmetric positive definite case one may use a block version of the strategy by~\cite{AjiJ84} in order to guarantee that the block incomplete
factorization does not break down\radd{ due to the presence of singular or ill-conditioned diagonal blocks}. 
\rsub{However
  i}{I}n the general case\radd{, on the other hand,} there exists no \rsub{analogy and even pivoting were required}{analogous strategy}.
\radd{Even in the symmetric positive definite case it was already observed in \cite{Ker78} that shifting the diagonal entry is already sufficient when there are not too many undesired pivots.}
Since our \rsub{block incomplete }{BI}LU approach does not use pivoting except inside
the diagonal blocks \radd{when employng LAPACK-based dense matrix kernels}, it may occasionally happen that diagonal blocks become
singular or \rsub{ill--conditioned}{ill-conditioned} in spite of having the system preprocessed
using maximum weight matching. To bypass this bottleneck (at least partially),
we perturb the diagonal blocks as follows\rsub{.}{:}
Let $\alpha=\max_{i,j} |a_{ij}|$ be  the maximum entry of $A$ (after scaling) in \rsub{modulus}{absolute value} and let $\tau$ and $\rho$ be some fixed
absolute and relative tolerance (in practice we use $\tau=10^{-2}$
and $\rho=10^{-1}$). Suppose that column $j$ of a diagonal block $D_{kk}$
consists of entries $d=(d_{ij})_{i=1,\dots,m}$. We denote their maximum
entry in \rsub{modulus}{absolute value} by $\delta_j$. If $d=0$ or if it turns out during
the $LU$ factorization that the block diagonal system
 is singular or ill--conditioned,
then we perturb the largest entry $d_{kj}$ in \rsub{modulus}{absolute value}  of $d$ 
by $d_{kj}^{(new)}=d_{kj}(1+\rho\beta_j)+\mathrm{sign}(d_{kj})\tau\alpha$. We give
preference to the diagonal entry instead of $d_{kj}$ (i.e.\radd{,} we choose $k=j$),
whenever $2|d_{jj}|\geqslant|d_{kj}|$. After that we proceed analogously
with respect to the rows of the diagonal block $D_{kk}$. By giving preference 
to the diagonal entries of $D_{kk}$ we reveal the original concept of maximum weight
matching. Moreover, \rsub{there is a slightly greater of hope that the system becomes
nonsingular or better condition with this kind of tie-break strategy}{this tie-breaking strategy might make the system nonsingular or of better condition} (e.g.\radd{,} consider
a matrix with entries of \radd{the} same order of magnitude that is rank-deficient). 
Perturbing the diagonal
blocks\radd{,} in general\radd{,} has to be \rsub{considered}{applied} with care and may easily introduce
severe numerical problems, but as long as the number of perturbations
is relatively small, this perturbation changes the factorization by a \rsub{low--rank}{low-rank} modification and the latter can usually be handled safely by Krylov subspace
methods. \radd{Alternatively to perturbing some diagonal blocks if necessary one could have restarted
  BILU applied to a shifted system which has been observed to be quite helpful \cite{Ben02}. However, in
  our comparisons we did not observe that BILU behaved better when using shifts.}

\subsection{Summarizing the Components of the Algorithm} \label{sect:sum}
After having explained the components that are combined to build up
the \rsub{block }{B}ILU we briefly summarize the main ingredients\rsub{.}{:}
\begin{enumerate}
\item Initially we apply maximum weight matching in order to
improve the diagonal dominance, i.e., $A\to \hat A=D_lAD_r\Pi$ (see \rsub{S}{s}ection~\ref{sect:mc64}).
\item Apply the \rsub{cosine--based}{cosine-based} blocking approach to $\hat A$ as described in
\rsub{S}{s}ection~\ref{sect:cosine}. This way we obtain from $\hat A$ 
a permuted matrix
$\tilde A = Q^T\hat A Q$.
\item Next reorder the compressed graph of $\tilde A$. Here the compressed
graph refers to the matrix $B$, where any diagonal block of $A$ according
to the \rsub{cosine--based}{cosine-based} blocking strategy is replaced by a scalar whenever
there is at least one nonzero entry inside this block.
We use nested dissection~\cite{karypis:98,LasK13} 
for reordering $B$ and we expand the permutation
afterwards in order to preserve the block structure of $\tilde A$.
From $\tilde A$ the next reordered matrix we obtain 
is $\check{A} = P^T \tilde A P$.
\item Given $\check A$, we simulate the behavio\rdel{u}r of our \rsub{block }{B}ILU
using the simplified $ILU(1,\tau)$ method from \rsub{S}{s}ection~\ref{sect:ilu1t}.
This simulation does not change $\check A$ anymore but it provides an
initial block structure prior to starting the \rsub{block }{B}ILU computation
\item Based on $\check A$ and its block structure, compute the \rsub{Crout--type}{Crout-type}
\rsub{block }{B}ILU according to drop tolerance $\tau$.
\item While computing the \rsub{block }{B}ILU, attempt to aggregate blocks progressively
in order to build larger blocks on the fly.
\end{enumerate}

\radd{Summing up all components we eventually end up with an approximate factorization
  $A\approx D_l^{-1}P_lLD^{-1}UP_r^TD_r^{-1}$ which will be used as preconditioner for Krylov
  subspace methods. Here, $D_l,D_r$ refer to the diagonal scaling matrices from \nref{eqn:mwm}, $P_l=QP$
and $P_r=\Pi QP$ are the permutation matrices collected from  \nref{eqn:mwm}, \nref{eqn:cos} and \nref{eqn:ordering} and $LD^{-1}U$ is the core BILU.}
It should be clear\rdel{,} that\radd{,} depending on the application\radd{,} one certainly
may skip one of these steps. E.g.\radd{,} maximum weight matching is\radd{,} in general\radd{,}
very beneficial as part of a \rsub{black--box}{black-box} approach (see\radd{,} e.g.\radd{,}~\cite{benzi:2000:phi})\rsub{,}{;}
however\radd{,} for some specific applications one might \radd{want} to avoid it 
because of its nonsymmetric permutation which is not always helpful.
Similarly, nested dissection is chosen just as one fill-reducing
ordering\radd{,}{;} other orderings such as \rsub{AMD}{approximate minimum degree (AMD)}~\cite{AmeDD96} could have been used as 
well. Also, e.g., the \rsub{cosine--based}{cosine-based} approach may not always pay off if the 
pattern of the original matrix does not have enough inherent block
structures. \rdel{We have chosen this set up in order to make the
approach at least half way a black box approach and to discuss the novelty
of the components as part of the complete block factorization approach.}
\radd{We have included this preprocessing procedure in the experiments for two reasons: first, to make the approach halfway a black-box approach, since the cosine-based might fail to provide improvement for unstructured problems; second, to discuss the novelty of the components as part of the complete block factorization approach.}


\section{Numerical Experiments} \label{sect:exp}
For the numerical experiments we select 100 (cf. Appendix~\ref{sect:app}) 
\rsub{large--scale}{large-scale} nonsymmetric
real nonsingular sparse matrices from the {SuiteSparse Matrix Collection}
(see Example~\ref{exm:venkat50-1}), each of them having a size of at least
$n\geqslant 50000$. Furthermore we use the same hardware configuration as
in Example~\ref{exm:venkat50-3}, which consists of a
single node with 1 TB main memory and $4$ Intel Xeon E7-4880 v2 @ 2.5 GHz 
processors
each of them having $15$ cores on a socket leading to $60$ cores in total.
As numerical methods we use the scalar Crout--type ILU as implemented
as binary code in \radd{the} MATLAB\rdel{'s} \texttt{ilu} (referred to as \rdel{``}ILUC\rdel{''}) and \radd{the} MATLAB\rdel{'s}
GMRES~\cite{SaaS86} implementation with a restart length of $30$ and relative residual
threshold $10^{-6}$. Our own variants of the \rdel{block ILU (termed ``}BILU\rdel{'')}
are implemented in C and use GMRES(30) as iterative solver as well.
In order
to distinguish between the single blocking strategies we add \rsub{in}{to} our results
suffixes such as ``-{}-p'' or ``cip'' in order to illustrate, which and
how many
of the three blocking strategies ``cosine'' (c), ``$ILU(1,\tau)$'' (i)\radd{,}
and ``progressive aggregation'' (p) are used in combination with BILU.
\radd{Notice that BILU(-~-~-) reduces to a scalar ILU.}
All matrices are preprocessed with maximum weight matching MC64~\cite{duko:99a} 
and reordered
with nested dissection METIS~\cite{karypis:98,LasK13} 
(in \radd{the} case that the cosine blocking is used, METIS is applied to the compressed graph). We use drop tolerances $\tau=10^{-1},10^{-2},\dots,10^{-6}$ and finally select the fastest ILU with respect to
this selection of $\tau$. It is clear for incomplete factorization methods
that their applicability is parameter\rsub{--}{ }dependent and if $\tau=10^{-1}$
is optimal for one system it may happen that $\tau=10^{-6}$ is required
for another system. To compensate \rsub{with}{for} the large variety of problems
we also state how often which choice of $\tau$ was selected.
Beside ILUC as one  benchmark we use PARDISO~\cite{sg:04-fgcs}
as another competitor, knowing
that over a large selection of matrices, direct solvers are typically known
to outperform iterative solvers. However, comparing with PARDISO
allows us to measure how far or how close the new \rsub{block }{B}ILU is regarding
an \rsub{up to date}{up-to-date} high performance sparse direct solver. Interestingly,
PARDISO uses maximum weight matchings and nested dissection initially as
well which makes the comparison even more appropriate.
Besides, we also compare the \rsub{block }{B}ILU with UMFPACK as implemented in MATLAB
and with SuperILU~\cite{LiS11} using a similar set up as for BILU.

  \subsection{Results}
  In order to evaluate the quality of the different incomplete
factorization methods, PARDISO and UMFPACK\radd{,} for the large selection of 
test problems, we
use performance profiles as a tool for benchmarking and for comparing
the algorithms. These\rdel{s} profiles were first proposed
in~\cite{dolan:2002} for benchmarking optimization software and
subsequently became the standard evaluation tool in the linear solver
and optimization community~\cite{scottgould:2004}. The profiles are
generated by running the set of methods $\mathcal{M}$ (eight variants of BILU,
ILUC, SuperILU, UMFPACK\radd{,} and PARDISO) on our
set of sparse matrices $\mathcal{S}$ and recording information of
interest, e.g.\radd{,} time for the solution operation for a required 
drop tolerance $\tau$ and memory consumption. 
Let us assume that a method $ m \in
\mathcal{M}$ reports a statistic $t_{ms} \ge 0$ for a matrix $s \in \mathcal{S}$ and that a smaller statistic
$t_{ms}$ indicates a better solution strategy. We can further define
$\tilde{t}_{s} = \min\{\;t_{ms}, m \in \mathcal{M}\; \}$, which
represents the best statistic for a given  matrix $m$.  Then for
$\alpha \ge 0$ and each $m \in \mathcal{M}$ and $s \in \mathcal{S}$ we
define
\begin{equation}
      k ( t_{ms}, \tilde{t}_{s}, \alpha) = \left \{ 
          \begin{aligned}
              1 \quad &\mbox{if} \; t_{ms} \le \alpha \, \tilde{t}_{s} \radd{\,,} \\
              0 \quad &\mbox{otherwise} \,.
          \end{aligned}
          \right. 
\end{equation} 
The performance profile $ p_m({\alpha}) $ of the method $m$ is then
defined by
\begin{equation}
p_m({\alpha}) = \frac{\sum_{s \in \mathcal{S}} \, {k ( t_{ms}, \tilde{t}_{s}, \alpha)}}{  | \mathcal{S} |  }.
\label{performance-profile}
\end{equation} 
Thus, the values of $p_m(\alpha)$ indicate the
fraction of all examples which can be solved within $\alpha$ times
the best strategy\radd{,} e.g.  $p_m(1)$ gives the fraction of which solution 
method $m$ is the most effective method and
 $\lim_{\alpha \rightarrow  \infty} $ 
indicates the fraction for which the algorithm succeeded.

\begin{figure}
 \begin{center}
\includegraphics[width=0.55\textwidth,height=0.4\textwidth]{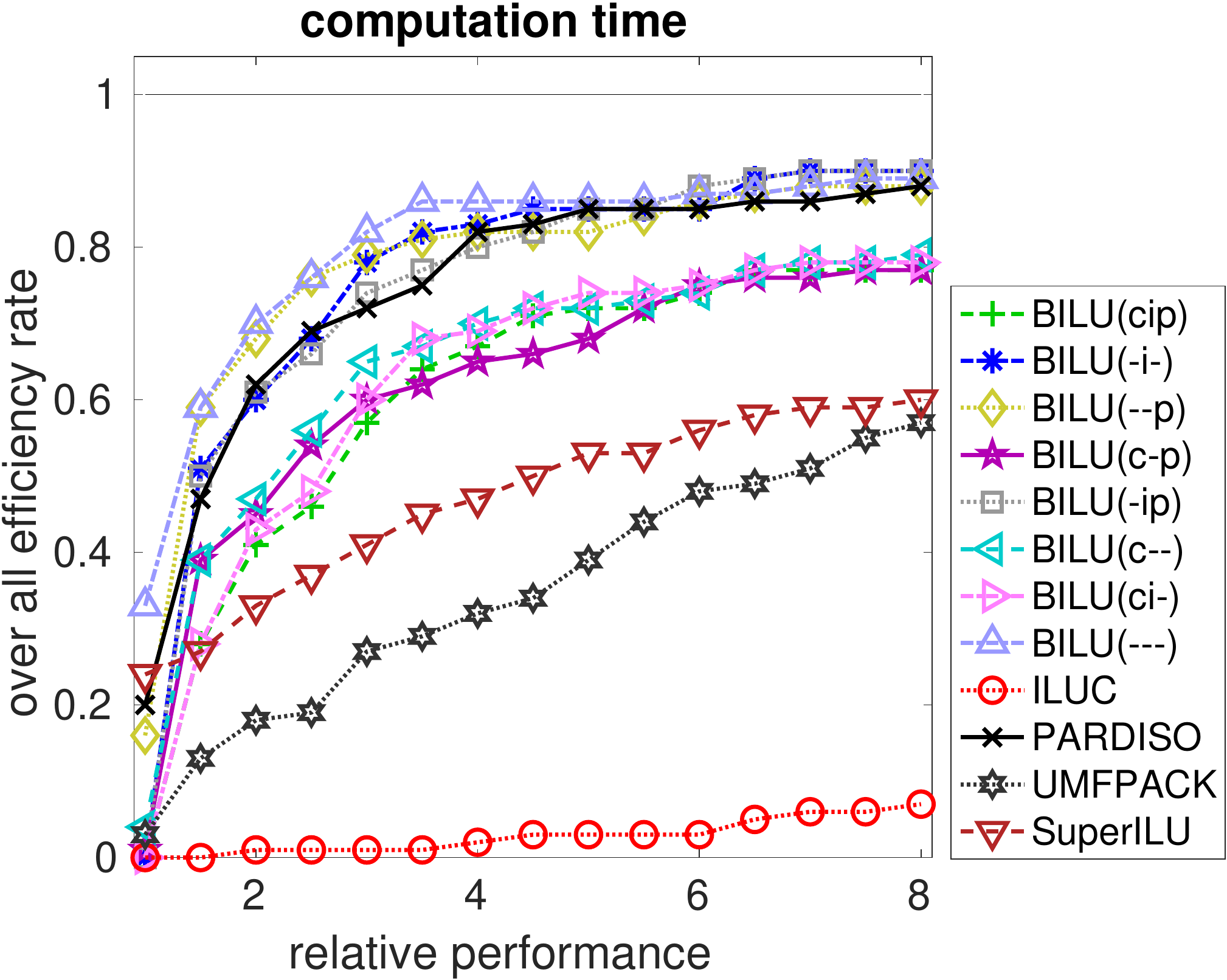} 
 \end{center}
\caption{Performance profile with respect to the best computation time.}
\label{fig:overall-time}
\end{figure}

\begin{figure}
 \begin{center}
\includegraphics[width=0.55\textwidth,height=0.4\textwidth]{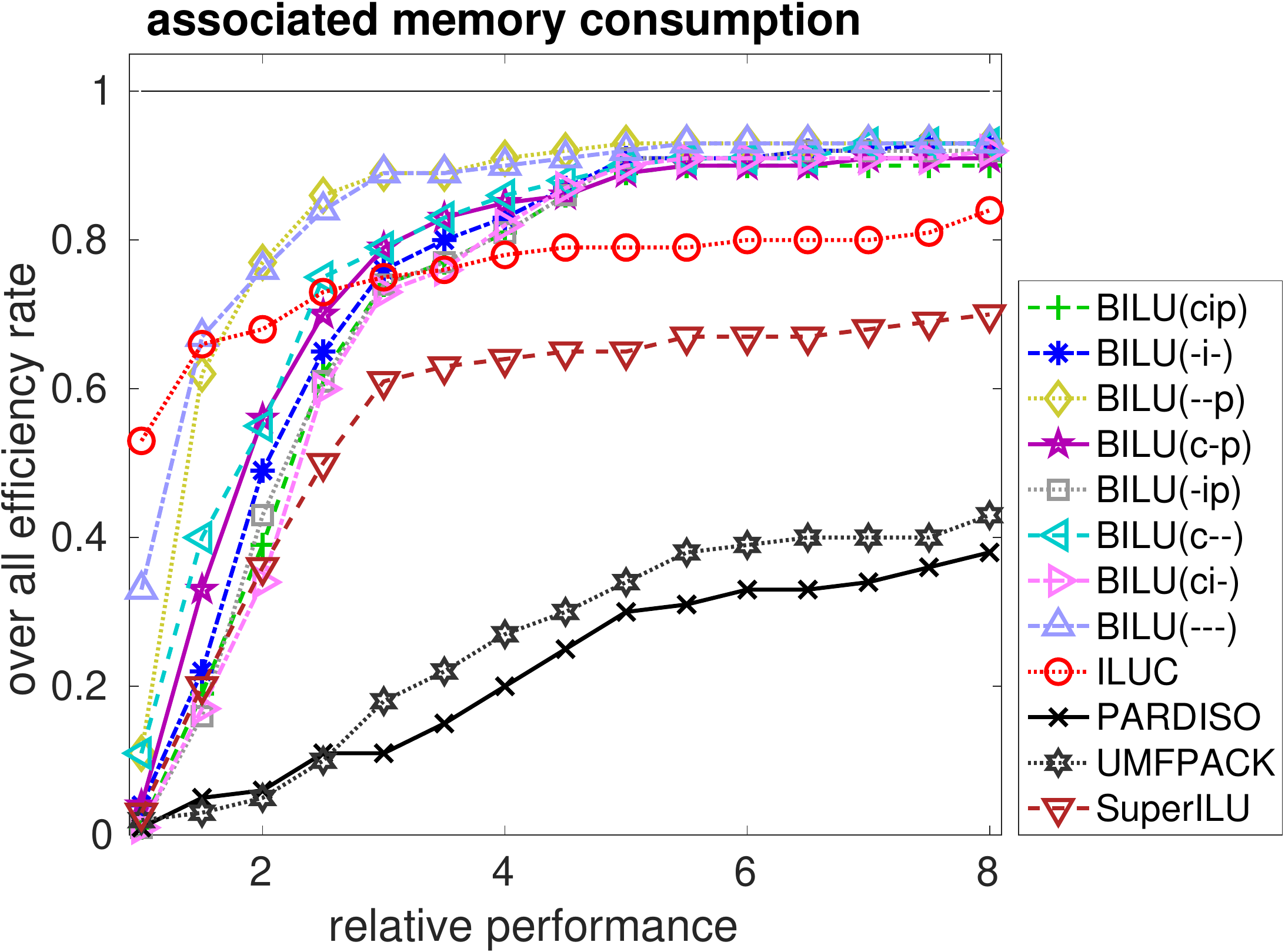} 
 \end{center}
\caption{Performance profile of the memory consumption associated with the best computation time.}
\label{fig:overall-mem}
\end{figure}

\begin{figure}
 \begin{center}
\includegraphics[width=0.55\textwidth,height=0.4\textwidth]{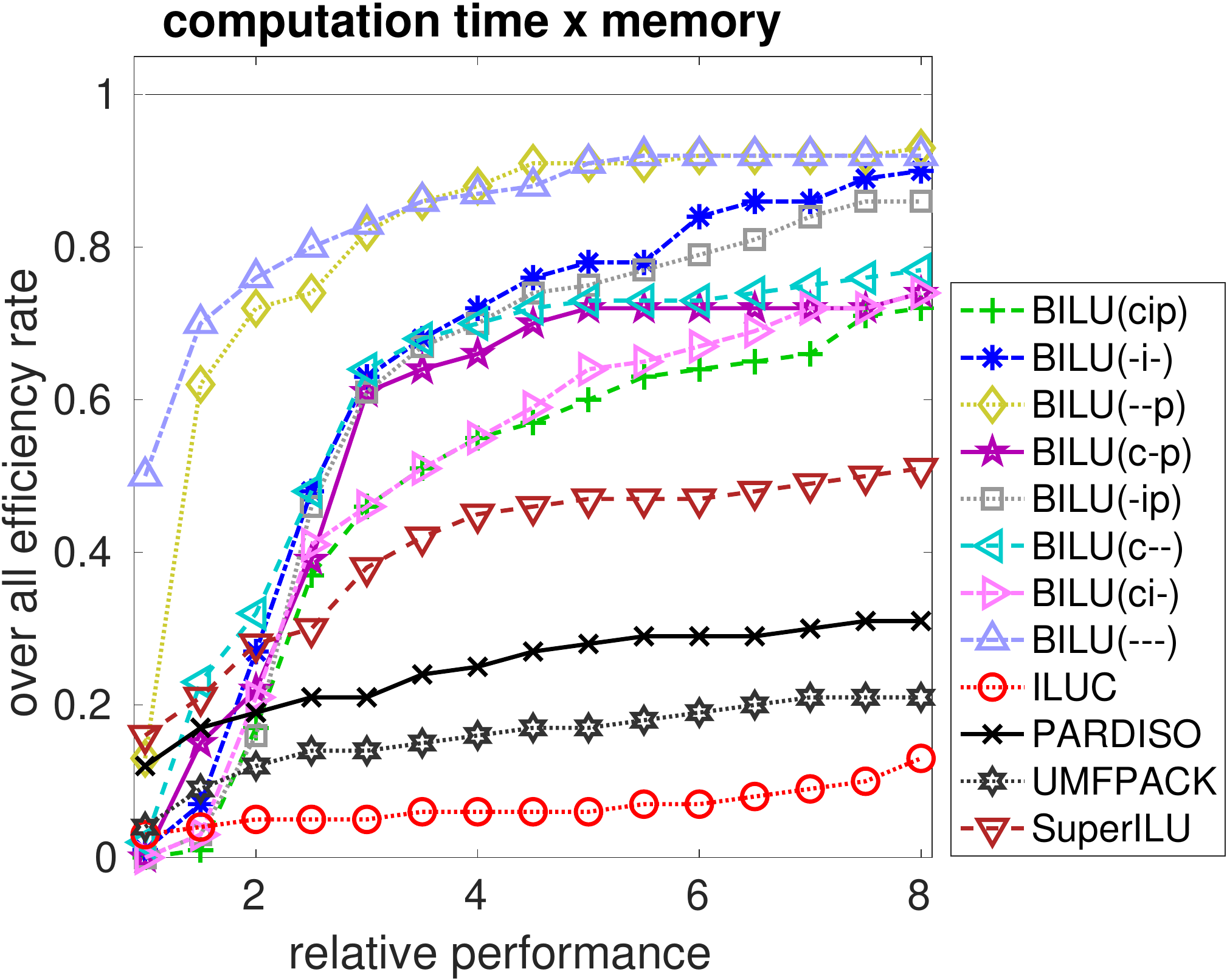} 
 \end{center}
\caption{Performance profile with respect to the best computation time $\times$ memory consumption.}
\label{fig:overall-timexmemory}
\end{figure}

\rdel{To report this statistics, we first display the best computation time in Figure~\ref{fig:overall-time}.}
\radd{To report these statistics, we first display the best computation time in Figure~\ref{fig:overall-time}.}
As we can easily see \rdel{from Figure~\ref{fig:overall-time}}, the \rsub{block }{B}ILU methods outperform the scalar
ILUC drastically. One has to be aware that \radd{the} block factorization method consumes more memory. In order
to demonstrate that the additional amount of memory is usually still acceptable, we display for
the methods from  Figure~\ref{fig:overall-time} the associated memory consumption as performance profile
in Figure~\ref{fig:overall-mem}.
As one would expect, Figure~\ref{fig:overall-mem} shows that the scalar factorization\rsub{ (here}{,} BILU(-~-~-)\rsub{)}{,} yields the smallest
amount of memory, but the variants of BILU using various blockings are most of the time within a close range of \radd{the} scalar version.
The use of approximate factorization methods as \radd{an} alternative to direct factorization methods is only partially
justified by their \rsub{fewer}{smaller} memory consumption. For many problems, as \rdel{a black--box}{black-box} solvers\radd{,} direct methods
are more reliable but occasionally too slow or too \rdel{much} memory consuming. A natural alternative statistics
is based on weighting memory and time appropriately by defining the best performance and the
product of time and memory~\cite{GeoGS12}. This performance profile is revealed in Figure~\ref{fig:overall-timexmemory}
showing that with respect to both aspects, time and memory, the BILU variants are apparently extremely
attractive.

\begin{figure}
 \begin{center}
\includegraphics[width=0.55\textwidth,height=0.4\textwidth]{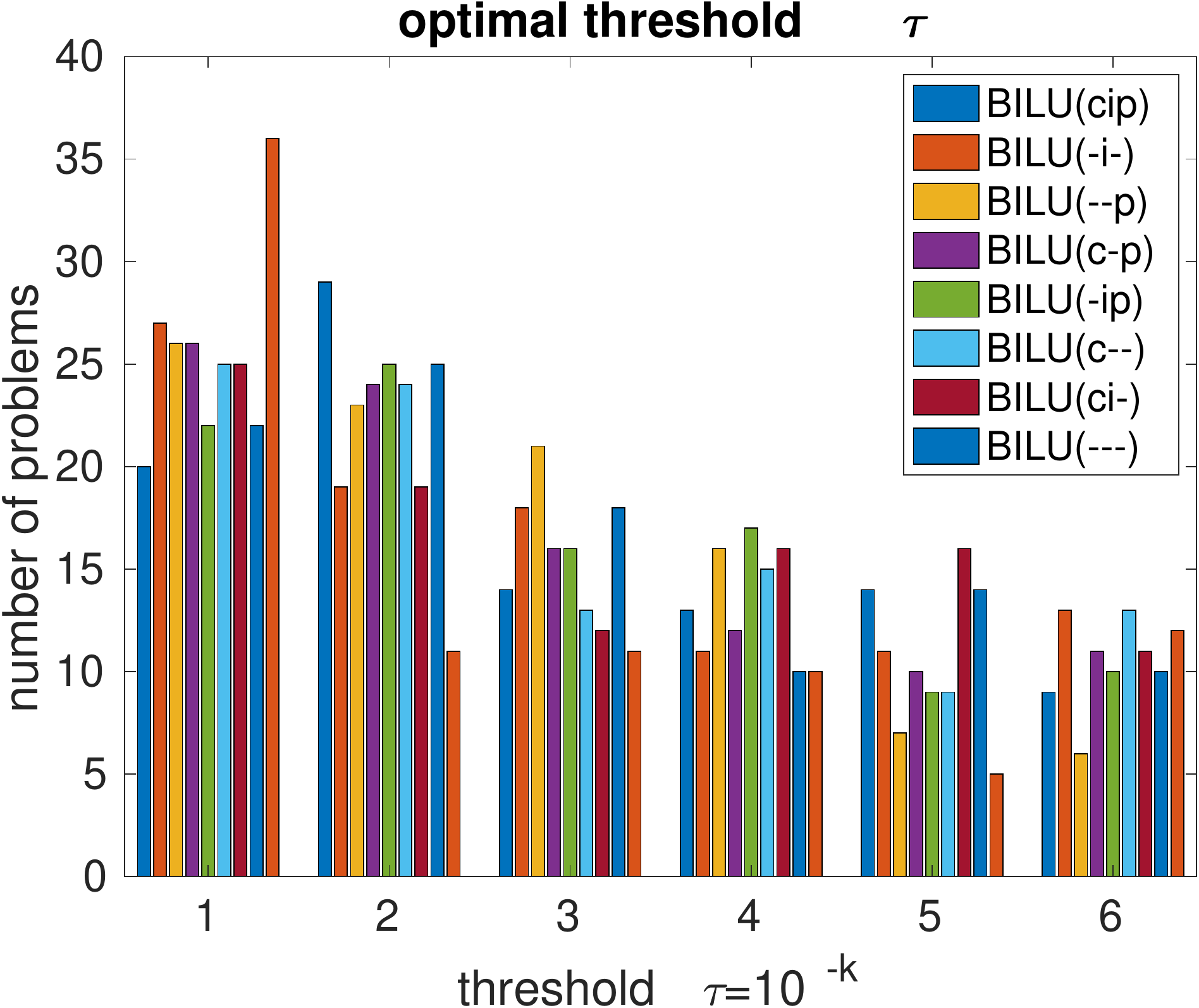} 
 \end{center}
\caption{Selection of drop tolerance $\tau$.}
\label{fig:overall-tau}
\end{figure}

We like to point out that small drop tolerances are rarely chosen which is in line
with the observations in~\cite{benzi:2000:phi}. This is illustrated in 
Figure~\ref{fig:overall-tau}.

\begin{figure}
 \begin{center}
\includegraphics[width=0.55\textwidth,height=0.4\textwidth]{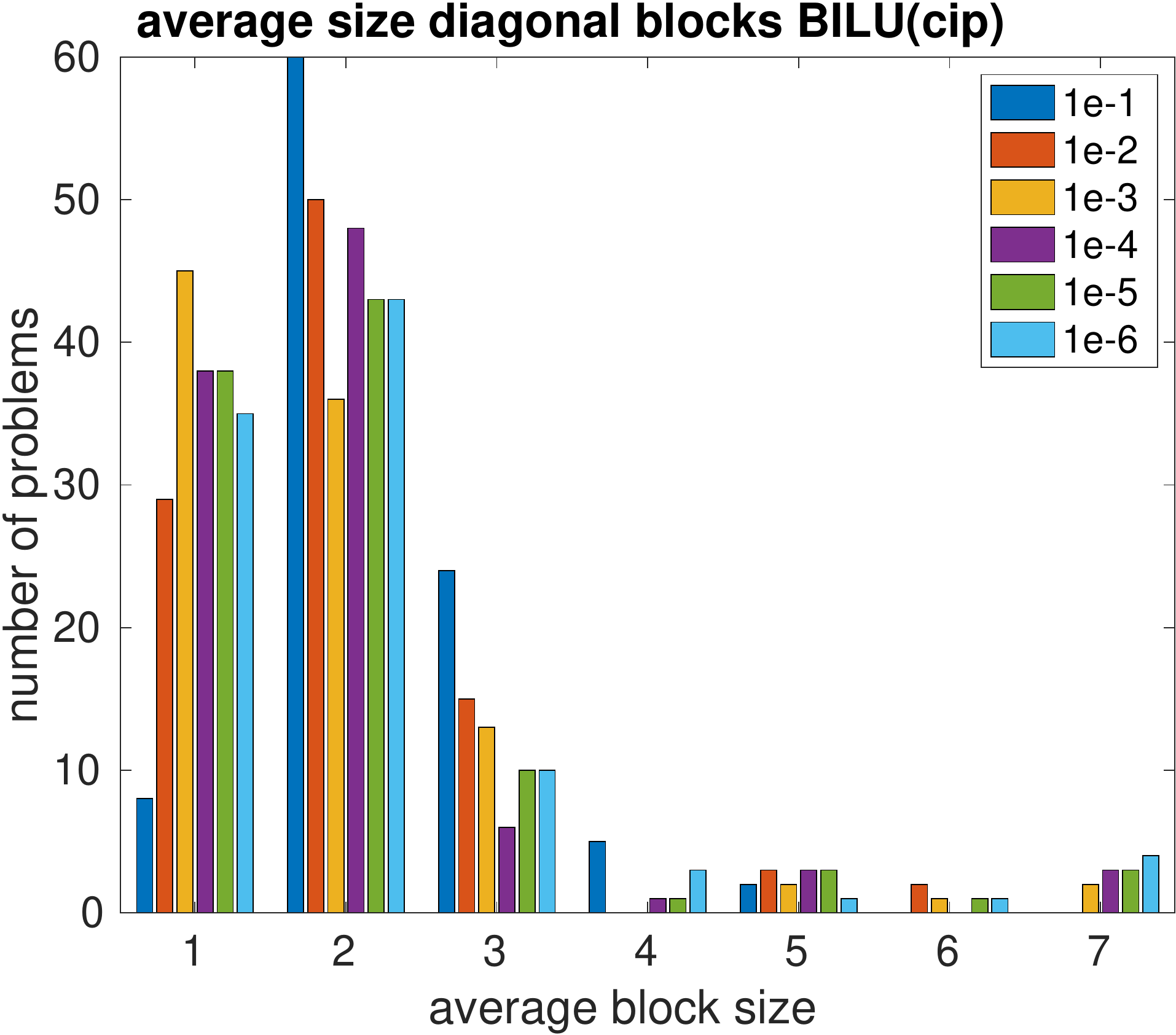} 
 \end{center}
\caption{Average block size for BILU(cip) with respect to the selected drop tolerance $\tau$.}
\label{fig:average-blocksize}
\end{figure}

Finally\radd{,} we stress that the large selection of application problem\radd{s} has le\rdel{a}d \rsub{to block--structured factorizations with many almost scalar structures, i.e., diagonal blocks of size $1$.}{the algorithm to select small sizes for the diagonal blocks, typically 1 and 2 (this analysis for BILU(cip) is reported in 
Figure~\ref{fig:average-blocksize}). While this is not true for structured problems, this is the average result when considering datasets of heterogeneous nature.}
\rdel{This is highlighted for BILU(cip) in 
Figure~\ref{fig:average-blocksize}.}
Having this almost ``scalar'' structure in mind, the block-structured approach is still very close
to the scalar version even in the frequent case when the factorization is relatively sparse and
nontrivial block structures occur rarely. This makes the \rsub{block--structured}{block-structured} approach competitive
even on a large scale of problems for which it is not optimally designed.

\subsection{Performance on selected problems}
In this section we consider \rsub{6}{six} real symmetric indefinite 
matrices (``af\_shell*'') which arise from industrial applications
in sheet metal forming. We compare the symmetric indefinite version 
(BILDL) of
our \rsub{block }{B}ILU which then becomes an incomplete block $LDL^T$ factorization.
Likewise, matching is replaced by a symmetrized approach as introduced
in~\cite{Duf2005}. \radd{Using~\cite{Duf2005}, a real diagonal matrix $D_{lr}$ and a permutation
matrix  $\Pi$ are computed such that
\begin{equation} \label{eqn:smwm}
\hat A=\Pi^TD_{lr}AD_{lr}\Pi
\end{equation}
and all entries of $\hat A$ satisfy $|\hat a_{ij}|\leqslant 1$. Moreover, in practice $\hat A$ will have
many diagonal blocks of size either $1\times 1$ such that $|\hat a_{ii}|=1$ or 
of size $2\times 2$ such that  $|\hat a_{i,i+1}|=|\hat a_{i+1,i}|=1$. For details we refer
to~\cite{Duf2005}. The cosine--based compression technique is then applied to the compressed
companion matrix, where the potential $2\times2$ pivots are merged. After that, compressing
the additional blocks from the cosine algorithm, a symmetric reordering is applied to the compressed
graph. The $ILU(1,\tau)$ is modified to deal with $1\times1$ and $2\times2$ pivots whatever is locally more appropriate.
This yields the symmetrically preprocessed block--structured matrix $\check A=P_{lr}^TD_{lr}AD_{lr}P_{lr}$,
where the permutation matrix $P_{lr}$ refers to the overall permutation. $\check A$ is approximately
factorized as $LDL^{T}$ using
the underlying block structure and a similar symmetrized perturbation strategy as in Section \ref{sect:perturbation}
is used whenever the diagonal blocks are ill-conditioned.}

\begin{figure}
 \begin{center}
\includegraphics[width=0.55\textwidth,height=0.4\textwidth]{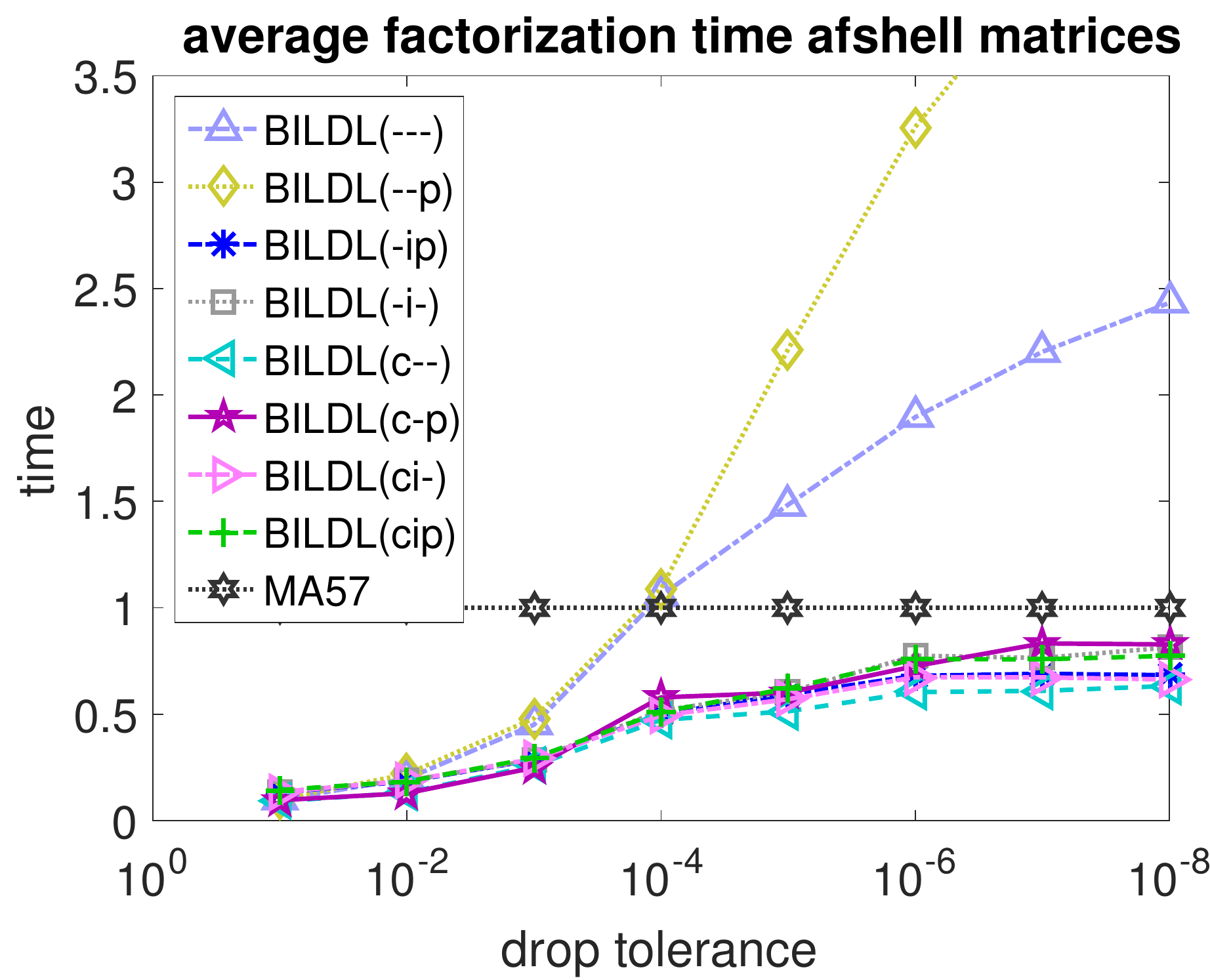} 
 \end{center}
\caption{Computation time of (block) incomplete $LDL^T$ factorizations compared with the symmetric indefinite direct solver MA57.}
\label{fig:performance-af-shell-bildl}
\end{figure}

We compare the \rsub{block--structured}{block-structured} incomplete
factorization approach with the direct solver MA57 as implemented in MATLAB.
Initially we compare the computation time for factorizing the matrix 
for the symmetric indefinite variants of BILU depending on the drop tolerances
with the computation time as required by MA57.
The results in 
Figure~\ref{fig:performance-af-shell-bildl} clearly demonstrate that
the scalar approaches are far out of competition whereas the block-structured
approach even remains \rsub{fast than}{competitive with} MA57 for relatively small drop tolerances
showing the other face of the \rsub{block--structured}{block-structured} approach, namely
turning more and more into a \rsub{high--performance}{high-performance} direct solver.
The computation time for each matrix is normalized by the smallest computation
time of BILDL and then averaged over the six sample matrices.

\begin{figure}
 \begin{center}
\includegraphics[width=0.55\textwidth,height=0.4\textwidth]{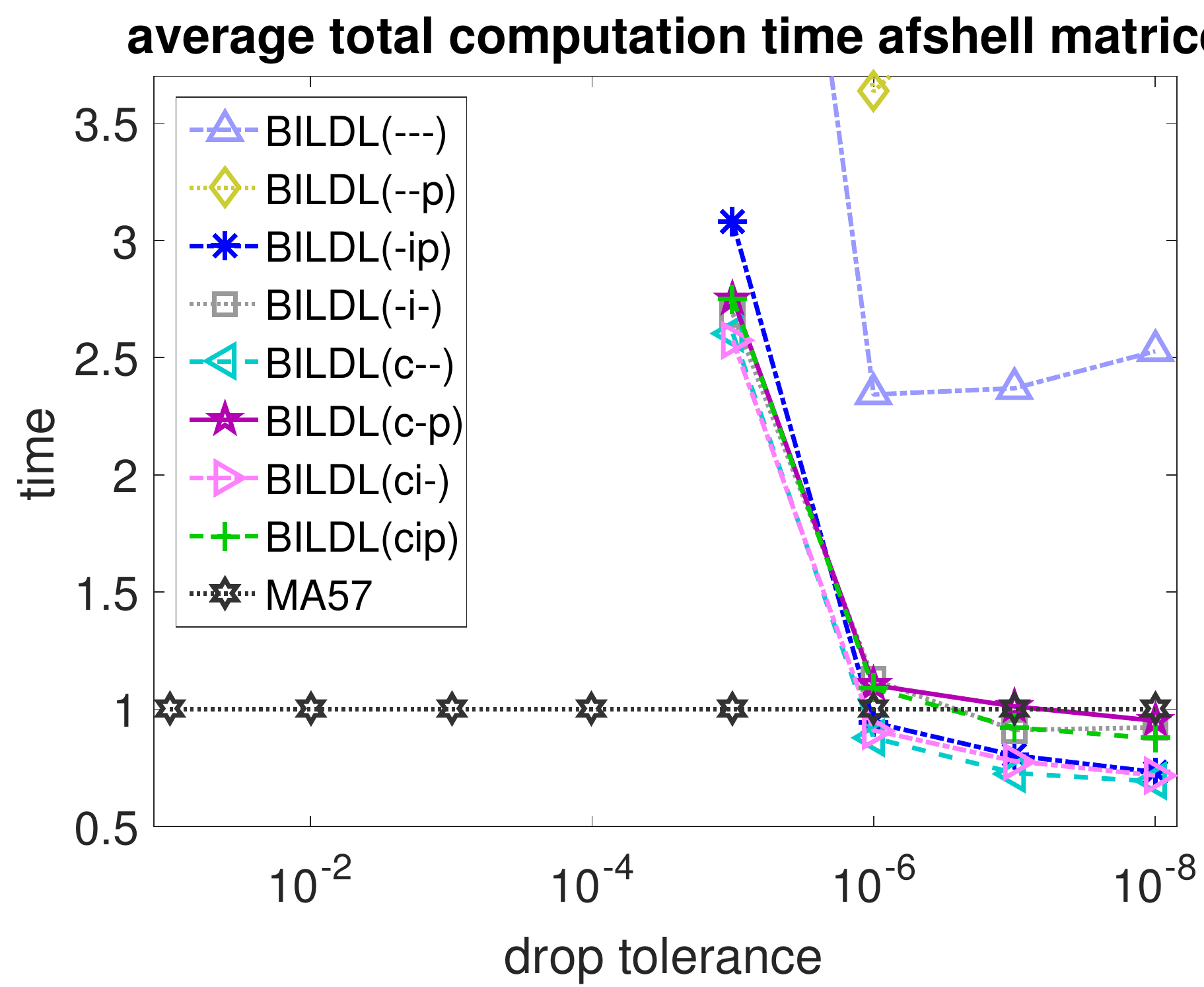} 
 \end{center}
\caption{Total computation time of (block) incomplete $LDL^T$ including SQMR compared with the symmetric indefinite direct solver MA57.}
\label{fig:performance-af-shell-bildl-total}
\end{figure}

Obviously, for preconditioning methods one has to incorporate the computation
time for the associated iterative solver, in our case we have chosen
the simplified QMR~\cite{FreJ97,FreN95} as iterative solver. This certainly
changes the situation since $\tau\leqslant10^{-5}$ was required in order
obtain convergence (we use the backward error and a tolerance of $10^{-6}$)\rsub{,}{;}
cf. Figure\rsub{.}{~}\ref{fig:performance-af-shell-bildl-total}

\begin{figure}
 \begin{center}
\includegraphics[width=0.55\textwidth,height=0.4\textwidth]{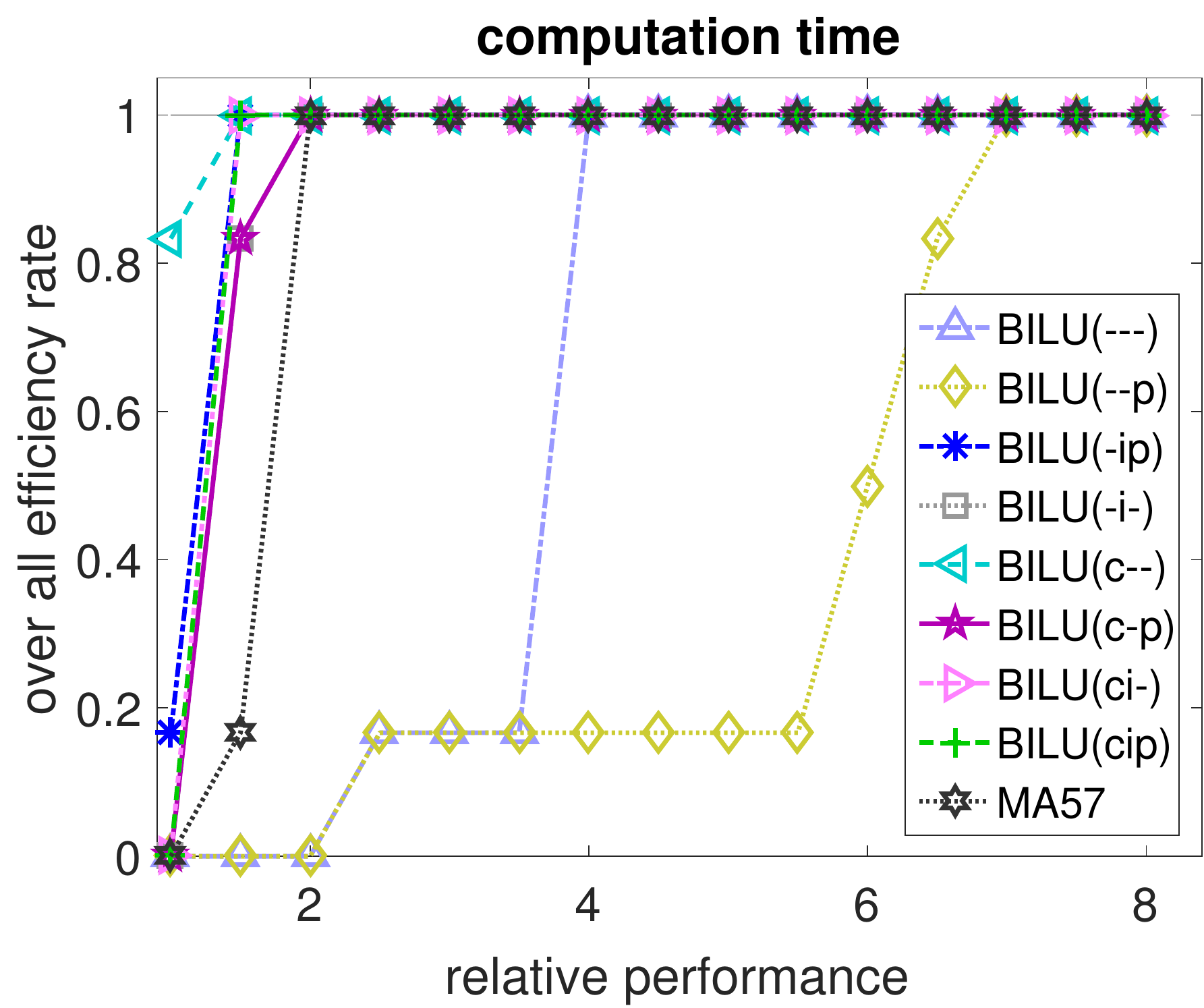} 
 \end{center}
\caption{Performance profile with respect to the best computation time, af\_shell matrices.}
\label{fig:performance-af-shell}
\end{figure}

\begin{figure}
 \begin{center}
\includegraphics[width=0.55\textwidth,height=0.4\textwidth]{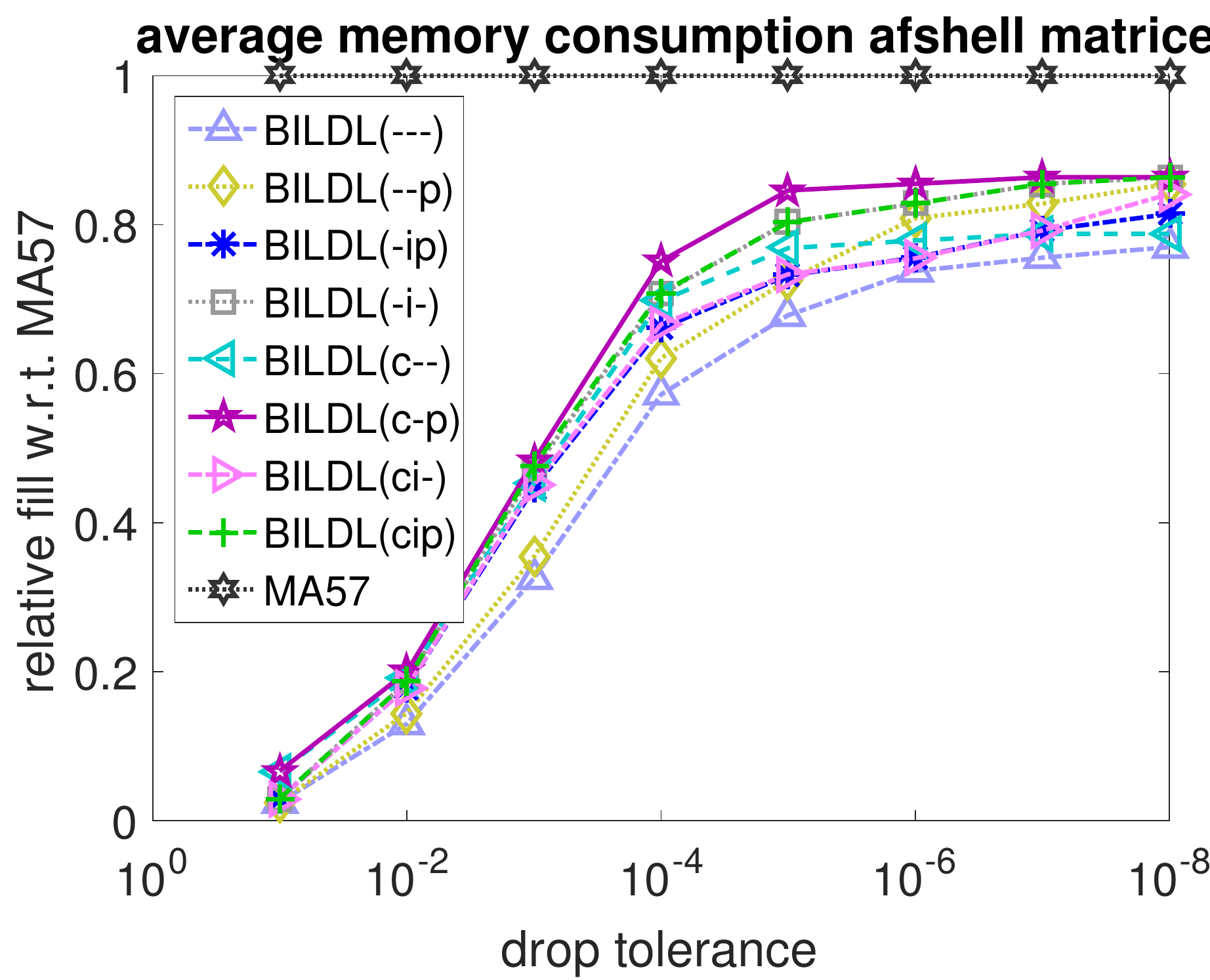} 
 \end{center}
\caption{Memory consumption block incomplete $LDL^T$ relative to MA57 (normalized to $1.0$).}
\label{fig:performance-mem-af-shell}
\end{figure}

In order to better display the total performance we draw a performance profile
(\ref{performance-profile})
in analogy to the previous section\rsub{,}{;} see Figure~\ref{fig:performance-af-shell}.
The performance profile clearly underlines the strength of the \rsub{block--structured}{block-structured}
approach even in comparison with a \rsub{high--performance}{high-performance} direct solver, whereas
the scalar version suffers from the large amount of \rsub{fill--in}{fill-in}.
This fill is illustrated in Figure~\ref{fig:performance-mem-af-shell} which demonstrates that the \rsub{block--structured}{block-structured}
ILU consumes memory close to the amount that is required by MA57, at least
for smaller drop tolerances.

\section{Concluding remarks}
We have demonstrated that using blocking strategies we are able to 
create a high performance incomplete \rsub{block }{B}ILU that is able to outperform
standard \rsub{incomplete }{$I$}$LU$ factorization by orders of magnitude on modern
computer architectures. Beside
the blocking strategies\radd{,} the use \radd{of} dense matrix kernels is the major reason
for its dramatic success \radd{in} closing the gap between ILUs and \rsub{up to date}{up-to-date} sparse
direct solvers. Beyond the scope of this paper is the integration of
BILU\radd{\footnote{JANUS BLOCK ILU available at \url{https://bilu.tu-bs.de}}}
 as template inside multilevel factorization methods. We plan to 
investigate this topic in the near future.


\bibliographystyle{abbrv}
\bibliography{bilu}

\section{Appendix}\label{sect:app}

List of matrices from the SuiteSparse Matrix Collection
used for the numerical experiments.

\begin{longtable}{lrrc@{~~}lrr}
name & size $n$ & $\frac{nz(A)}n$        &&   name & size $n$ & $\frac{nz(A)}n$\\
\cline{1-3}\cline{5-7}
\endhead
2D\_54019\_highK & $  54019$ & $9.0$&&        lhr71c & $  70304$ & $    21.7$ \\
3D\_51448\_3D & $  51448$ & $  10.4$&&        lung2 & $ 109460$ & $     4.5$  \\       
ASIC\_100k  & $  99340$ & $     9.5$&&        majorbasis &$ 160000$ & $10.9$  \\
ASIC\_100ks & $  99190$ & $     5.8$&&        mark3jac120 & $  54929$ & $     5.9$\\    
ASIC\_320k  & $ 321821$ & $     6.0$&&        mark3jac120sc & $  54929$ & $   5.9$\\  
ASIC\_320ks & $ 321671$ & $     4.1$&&        mark3jac140 & $  64089$ & $     5.9$\\   
ASIC\_680k  & $ 682862$ & $     3.9$&&        mark3jac140sc & $  64089$ & $   5.9$\\  
ASIC\_680ks & $ 682712$ & $     2.5$&&        matrix\_9 & $ 103430$ & $    11.7$\\ 
atmosmodd   & $1270432$ & $     6.9$&&       matrix-new\_3 & $ 125329$ & $   7.1$\\         
atmosmodj   & $1270432$ & $     6.9$&&       memchip & $2707524$ & $     4.9$\\    
atmosmodl   & $1489752$ & $     6.9$&&      ohne2  & $181343$ & $37.9$\\       
barrier2-1  & $ 113076$ & $    18.8$&&         para-4 & $153226 $ & $19.1 $ \\        
barrier2-2  & $ 113076$ & $    18.8$&&           para-5  & $ 155924 $ & $13.4 $ \\      
barrier2-3  & $ 113076$ & $    18.8$&&        para-6  & $155924 $ & $13.4 $ \\ 
barrier2-4  & $ 113076$ & $    18.8$&&              para-7   & $155924 $ & $13.4 $ \\    
barrier2-9  & $ 115625$ & $    18.7$&&          para-8 & $ 155924$ & $    13.4$\\          
barrier2-10 & $ 115625$ & $    18.7$&&         para-9 & $ 155924$ & $    13.4$\\            
barrier2-11 & $ 115625$ & $    18.7$&&         para-10 & $155924$  &  $13.4$ \\            
barrier2-12 & $ 115625$ & $    18.7$&&         poisson3Db & $  85623$ & $    27.7$\\       
Baumann     & $ 112211$ & $     6.7$&&       Raj1    & $263743 $ & $4.9 $ \\              
bayer01      & $  57735$ & $     4.8$&&      rajat16 & $  94294$ & $     5.1$\\  
bcircuit     & $  68902$ & $     5.5$&&          rajat17 & $  94294$ & $     5.1$\\        
cage12       & $ 130228$ & $    15.6$&&          rajat18 & $  94294$ & $     5.1$\\                   
cage13       & $ 445315$ & $    16.8$&&              rajat20 & $  86916$ & $     7.0$\\              
cage14       & $1505785$ & $    18.0$&&           rajat21 & $ 411676$ & $     4.6$\\        
cage15       & $5154859$ & $    19.2$&&              rajat23 & $ 110355$ & $     5.0$\\      
Chebyshev4   & $  68121$ & $    78.9$&&              rajat24 & $ 358172$ & $     5.4$\\  
circuit\_4   & $  80209$ & $     3.8$&&            rajat25 & $  87190$ & $     7.0$\\ 
circuit5M\_dc& $3523317$ & $4.2$&&                rajat28 & $  87190$ & $     7.0$\\  
circuit5M    & $5558326$ & $    10.7$&&       rajat29 & $ 643994$ & $     5.8$\\     
crashbasis   & $ 160000$ & $    10.9$&&      rajat30 & $ 643994$ & $     9.6$\\
dc1          & $ 116835$ & $     6.6$&&       rajat31 & $4690002$ & $     4.3$\\                    
dc2          & $ 116835$ & $     6.6$&&          scircuit  & $170998 $ & $5.6 $ \\               
dc3          & $ 116835$ & $     6.6$&&                 shyy161 & $  76480$ & $     4.3$\\       
ecl32        & $  51993$ & $     7.3$&&               stomach     & $ 213360$ & $14.2 $ \\     
  epb3         & $  84617$ & $     5.5$ &&             tmt\_unsym  & $917825 $ & $5.0$ \\
FEM\_3D\_thermal2 &$147900$& $23.6$  &&        torso1 & $ 116158$ & $    73.3$\\                    
Freescale1   & $3428755$ & $     5.0$&&                torso2 & $ 115967$ & $     8.9$\\        
FullChip     & $2987012$ & $     8.9$&&      torso3 & $ 259156$ & $    17.1$\\                 
g7jac180     & $  53370$ & $    12.0$&&           trans4 & $ 116835$ & $     6.4$\\               
g7jac180sc   & $  53370$ & $    12.0$&&           trans5 & $ 116835$ & $     6.4$\\
g7jac200     & $  59310$ & $    12.1$&&       transient & $ 178866$ & $     5.4$\\             
g7jac200sc   & $  59310$ & $    12.1$&&        TSOPF\_RS\_b39\_c30 & $60098$ & $18.0$\\              
hcircuit     & $ 105676$ & $     4.9$&&       twotone & $ 120750$ & $    10.0$\\                  
hvdc2        & $189860$ &  $7.1$   &&          venkat01 & $  62424$ & $    27.5$\\                
ibm\_matrix\_2 & $  51448$ & $10.4$&&       venkat25 & $  62424$ & $    27.5$\\              
laminar\_duct3D& $  67173$ & $56.4$&&          venkat50 & $  62424$ & $    27.5$\\                
language       & $399130$ & $3.0$ &&                   water\_tank & $  60740$ & $    33.5$\\ 
largebasis     & $440020$ & $11.9$ &&          Wordnet3 & $  82670$ & $     1.6$\\                
lhr71          & $  70304$ & $    21.3$&&         xenon2     & $157464 $ & $24.6$ \\     

\end{longtable}

\end{document}